\documentclass{article}
\usepackage[utf8]{inputenc}

\usepackage{graphicx, amssymb, latexsym, amsfonts, amsmath, lscape, amscd, multirow, multicol,
  amsthm, color, epsfig, mathrsfs, tikz, enumerate, soul, subcaption, scrextend, todonotes, xstring, hyperref, comment}
\usetikzlibrary{calc}
\usepackage{array}
\usepackage{makecell}
\theoremstyle{definition}

\newtheorem{theorem}{Theorem}
\newtheorem{conjecture}[theorem]{Conjecture}

\newtheorem{corollary}[theorem]{Corollary}

\newtheorem{lemma}[theorem]{Lemma}
\newtheorem{proposition}[theorem]{Proposition}

\newtheorem{definition}[theorem]{Definition}
\newtheorem{observation}[theorem]{Observation}

\newtheorem{claim}{Claim}
\renewcommand*{\theclaim}{\Alph{claim}}

\theoremstyle{definition}

\newtheorem{subclaim}{Claim}[claim]

\newtheorem{subsubclaim}{Claim}[subclaim]

\newcommand{\smallo}{{\rm o}}

\newcommand{\even}{{\rm even}}
\newcommand{\odd}{{\rm odd}}

\newcommand{\barA}{{\overline{A}}}

\newcommand{\smallqed}{{\tiny ($\Box$)}}
\newcommand{\QED}{$\Box$}
\newcommand{\1}{\vspace{0.1cm}}
\newcommand{\2}{\vspace{0.2cm}}

\let\oldenumerate\enumerate
\renewcommand{\enumerate}{
  \oldenumerate
  \setlength{\itemsep}{0pt}
  \setlength{\parskip}{0pt}
  \setlength{\parsep}{0pt}
}

\newcommand\DELETE[1]{}

\newcommand{\ID}{\gamma^{{\rm ID}}}

\title{Identifying codes in triangle-free graphs of bounded maximum degree
\footnote{The first two authors were supported by the French government IDEX-ISITE initiative CAP 20-25 (ANR-16-IDEX-0001), the International Research Center "Innovation Transportation and Production Systems" of the I-SITE CAP 20-25, and the ANR project GRALMECO (ANR-21-CE48-0004). The research of Michael Henning was supported in part by the South African National Research Foundation (grant number 129265) and the University of Johannesburg.}}

\author{Dipayan Chakraborty\footnote{Corresponding author: dipayancha@gmail.com}~\footnote{\noindent Université Clermont Auvergne, CNRS, Clermont Auvergne INP, Mines Saint-Etienne, LIMOS, 63000 Clermont-Ferrand, France.}~\footnote{\noindent Department of Mathematics and Applied Mathematics, University of Johannesburg, South Africa.}
\and Florent Foucaud\footnotemark[3]
\and Michael A. Henning\footnotemark[4]
\and Tuomo Lehtil\"a\footnotemark[3]~\footnote{University of Helsinki, Department of Computer Science, Helsinki, Finland.}~\footnote{Helsinki Institute for Information Technology (HIIT), Espoo, Finland.}~\footnote{University of Turku, Department of Mathematics and Statistics, Turku, Finland.}~~\footnote{Research partially supported by Business Finland project 10769/31/2022 6GTNF and Academy of Finland grants 338797 and 358718.}
}

\date{}

\begin{document}

\maketitle

\begin{abstract}
An \emph{identifying code} of a closed-twin-free graph $G$ is a set $S$ of vertices of $G$ such that any two vertices in $G$ have a distinct, nonempty intersection between their closed neighborhood and $S$. It was conjectured that there exists a constant $c$ such that for every connected closed-twin-free graph $G$ of order $n$ and maximum degree $\Delta$, the graph $G$ admits an identifying code of size at most $\left( \frac{\Delta-1}{\Delta} \right) n+c$. In [D. Chakraborty, F. Foucaud, M. A. Henning, and T. Lehtil\"{a}. Identifying codes in graphs of given maximum degree: Characterizing trees. \emph{Discrete Mathematics}, 349(2):114826, 2026], we proved the conjecture for all trees. In this article, we show that the conjecture holds for all triangle-free graphs, with the same list of exceptional graphs needing $c>0$ as for trees: for $\Delta\ge 3$, $c=1/3$ suffices and there is only a set of 12 trees requiring $c>0$ for $\Delta=3$, and when $\Delta\ge 4$ this set is reduced to the $\Delta$-star only. Our proof is by induction, whose starting point is the above result for trees. Along the way, we prove a generalized version of Bondy's theorem on induced subsets [J. A. Bondy. Induced subsets. Journal of Combinatorial Theory, Series B, 1972] that we use as a tool in our proofs. We also use our main result for triangle-free graphs, to prove the upper bound $\left( \frac{\Delta-1}{\Delta} \right) n+1/\Delta+4t$ for graphs that can be made triangle-free by the removal of $t$ edges.
\end{abstract}

\section{Introduction}

In this article, we consider simple undirected loopless graphs. The \textit{open neighborhood} of vertex $u$ in graph $G$, $N_G(u)$, contains every vertex adjacent to vertex $u$ while the \textit{closed neighborhood} also contains vertex $u$ itself.
A set of vertices $D\subseteq V(G)$ is a \textit{dominating set} if every vertex outside of $D$ is adjacent to a vertex in $D$. A set $S\subseteq V(G)$ is a \textit{separating set} if every vertex $u\in V(G)$ has a unique closed neighborhood in set $S$, that is, intersection $N_G[u]\cap S$ is unique for every vertex. Furthermore, a set $C\subseteq V(G)$ is an \textit{identifying code} if it is a dominating set and a separating set \cite{karpovsky1998new}. See Figure~\ref{fig:trees} for examples of identifying codes. From the intuitive perspective, identifying codes allow us to locate or identify any vertex if we know which vertices in the set $S$ are in its closed neighborhood. It is natural to ask what is the minimum number of vertices in an identifying code of graph $G$. This value is called the \textit{identification number} of graph $G$ and it is denoted by $\ID(G)$.

This paper concentrates on proving Conjecture~\ref{conj_G Delta_UB} (stated below) on upper bounds for the minimum size of identifying codes of given maximum degree, for all triangle-free graphs. Previously, in~
\cite{PaperPart1}, we have proved the conjecture for trees. 

A thorough treatise on domination in graphs can be found in~\cite{HaHeHe-20,HaHeHe-21,HaHeHe-23}. Bounds on domination numbers for graphs with restrictions on their degree parameters are a natural and important line of research, see for example the influential result by Reed~\cite{reed_1996} that if $G$ is a connected cubic graph of order~$n$, then its domination number is at most $\frac{3}{8}n$. A detailed discussion on upper bounds on the domination number in terms of its order and degree parameters, as well as bounds with specific structural restrictions imposed, can be found in~\cite[Chapters 6, 7 and 10]{HaHeHe-23}.

Originally, identifying codes were motivated by fault-detection in multiprocessor networks \cite{karpovsky1998new}. Numerous other applications have been discovered, such as threat location in facilities using sensor networks~\cite{UTS04}, logical definability of graphs~\cite{PVV06} and canonical labeling of graphs for the graph isomorphism problem~\cite{B80}. Moreover, multiple related concepts have been introduced since the 1960s, such as \emph{separating systems}, \emph{test covers} and \textit{locating-dominating sets}, which have been independently discovered and studied, forming the general area of identification problems in graphs and other discrete structures. See for example~\cite{B72,MS85,R61,rall1984location}. An extensive internet bibliography containing over 500 articles around these topics can be found at \cite{lobstein2012watching}, while more information specifically about identifying codes can be found in the book chapter~\cite{lobstein2020locating}.

In this paper, we concentrate on connected identifiable triangle-free graphs. A graph is \textit{identifiable} if it does not contain any \textit{closed twins}, that is, vertices with the same closed neighborhoods. Moreover, a graph is \textit{triangle-free} if no three vertices in it form a cycle. Similarly as closed twins, we define \textit{open twins} as vertices with the same open neighborhoods. A graph without closed or open twins is called \textit{twin-free}. We also define similarly \textit{open-twin-free} and \textit{closed-twin-free} graphs. Twins are significant for separating sets and identifying codes; if two vertices are closed twins, then the graph does not admit any separating set and hence no identifying code. Moreover, if $t$ vertices have the same open neighborhood, then any separating set and hence also any identifying code, contains at least $t-1$ of them. Note that every connected triangle-free graph on at least three vertices is identifiable.

As our main result, we prove the following conjecture for all identifiable triangle-free graphs.

\begin{conjecture}[{\cite[Conjecture 1]{foucaud2012size}}] \label{conj_G Delta_UB}
There exists a constant $c$ such that for every connected identifiable graph $G$ on $n$ vertices and of maximum degree $\Delta\ge 2$,
\[
\ID(G) \le \left( \frac{\Delta - 1}{\Delta} \right) n + c.
\]
\end{conjecture}

\noindent 
In~\cite{PaperPart1}, we have proved Conjecture~\ref{conj_G Delta_UB} for trees and have determined the exact set of trees requiring a positive constant $c$ together with an exact value of $c$. In this paper, we use our previous results for trees from~\cite{PaperPart1} as a starting point for our proof of Conjecture~\ref{conj_G Delta_UB} for all triangle-free graphs.

It is known that if true, the conjecture would be tight, that is, some graphs of maximum degree $\Delta$ have identification number at least $\frac{\Delta-1}{\Delta}n$. For $\Delta=2$, the conjecture is tight for both paths and cycles with $c\le 3/2$ (see corollary~\ref{cor_paths & cycles}). For $\Delta=3$, the conjecture is tight for example for trees presented in Figure~\ref{fig:trees} and for a path whose every vertex we join to a $2$-path by a single edge (see~\cite{PaperPart1}). For any $\Delta>3$, the complete bipartite graph $K_{\Delta,\Delta}$ satisfies $\ID(K_{\Delta,\Delta})=2\Delta-2 = \left( \frac{\Delta-1}{\Delta} \right) n$ and hence gives tight examples with $c=0$. Furthermore, for any $\Delta>3$ and an unbounded value of $n$, there are trees with identification number $\ID(T)> \frac{(\Delta-1)n}{\Delta}-\frac{n}{\Delta^2}$ \cite{PaperPart1}. Furthermore, there exist connected graphs of any maximum degree~$\Delta\ge 3$ and arbitrarily large number $n$ of vertices, with identification number $\left( \frac{\Delta - 1}{\Delta} \right) n$, see~\cite{foucaud2012degree,sierpinski}.

The condition on the maximum degree $\Delta$ is a necessary part of the conjecture. Without it, there are graphs on $n$ vertices with identification number $n-1$~\cite{FGKNPV11,GM07}.

The best known general upper bounds for connected graphs with
a maximum degree $\Delta$ and number of vertices $n$, when $n$ is large enough, are of the form $n-\frac{n}{\Theta(\Delta^5)}$ \cite{FGKNPV11} which has been improved to $n-\frac{n}{103(\Delta+1)^3}$ in~\cite{foucaud2012degree} (for the sake of comparison, the conjectured bound can be rewritten as $n-\frac{n}{\Delta}+c$). When we consider graph classes instead of general graphs, some improvements on these bounds are known. If every pair of closed neighborhoods in a graph differ by at least two vertices, then the general bound has been improved to $n-\frac{n}{103\Delta}$ and to $n-\frac{n}{f(k)\Delta}$ for graphs of clique number $k$~\cite{foucaud2012degree}. For bipartite graphs, an upper bound of $n-\frac{n}{\Delta+9}$ has been proved \cite{foucaud2012size}. In a short conference proceedings paper \cite{IDBipartite}, we sketched a proof for Conjecture~\ref{conj_G Delta_UB} for bipartite graphs without twins of degree~2 or greater (we did not include that proof in the current article, since the triangle-free result we present here is stronger). Moreover, in~\cite{PaperPart1}, we gave the proof for all trees. Conjecture~\ref{conj_G Delta_UB} also holds for line graphs of graphs of average degree at least~5~\cite[Corollary 21]{line} as well as graphs which have girth at least $5$, minimum degree $2$ and maximum degree at least $4$ \cite{balbuena2015locating}. Furthermore, the conjecture holds in many cases for some graph products such as Cartesian and direct products~\cite{goddard2013id, junnila2019conjecture, rall2014identifying}. See also the book chapter~\cite{lobstein2020locating}, where Conjecture~\ref{conj_G Delta_UB} is presented.

Conjecture~\ref{conj_G Delta_UB} has been considered also for triangle-free graphs previously. In~\cite{foucaud2012size}, an upper bound of type  $n-\frac{n}{\Delta+\smallo(\Delta)}$ was presented for triangle-free graphs. When the triangle-free graph is also twin-free, this upper bound has been improved to $n-\frac{n}{3\Delta/(\ln\Delta-1)}$. Note that the latter result implies that Conjecture~\ref{conj_G Delta_UB} holds for triangle-free graphs without any open twins, whenever $\Delta\ge 55$ (because then, $3\Delta/(\ln\Delta-1)\leq\Delta$). Note that the graphs containing open twins seem to be the toughest cases regarding Conjecture~\ref{conj_G Delta_UB} (among triangle-free graphs). Indeed, we will see that every triangle-free graph requiring a positive constant $c$ with $\Delta\geq3$ contains open twins. Furthermore, in \cite{FT2022} it has been shown that every \emph{twin-free} bipartite graph $G$ on $n\geq5$ vertices satisfies the upper bound $\ID(G)\le \frac{2n}{3}$, while if we allow open twins, then there exist trees $T$ with arbitrarily large numbers of vertices such that $\ID(T) > \left( \frac{\Delta - 1}{\Delta} \right) n-\frac{n}{\Delta^2}$~\cite{PaperPart1}.

Our main result is to prove Conjecture~\ref{conj_G Delta_UB} (in a strong form) for all triangle-free graphs. To state it, we define, for every integer $\Delta\ge 3$, a set $\mathcal{F}_\Delta$ of exceptional graphs of maximum degree at most $\Delta$ (see Section~\ref{secExtremal} for greater detail). For $\Delta=3$, this set contains twelve trees (see Figure~\ref{fig:trees}), the cycles on 4 and 7 vertices, and the path on 4 vertices. For every integer $\Delta>3$, it contains exactly the $\Delta$-star $K_{1,\Delta}$.

\begin{theorem}\label{thm:trianglefree}
Let $\Delta\ge 3$ be an integer, and let $G$ be a connected triangle-free graph of order $n \ge 3$. If $G \in \mathcal{F}_{\Delta}$, then
\[
\ID(G) = \left( \frac{\Delta - 1}{\Delta}\right) n+\frac{1}{\Delta}.
\]
On the other hand, if $G \notin \mathcal{F}_{\Delta}$ has maximum degree $\Delta$, then
\[
\ID(G) \le \left( \frac{\Delta - 1}{\Delta}\right) n.
\]
\end{theorem}
Note that graphs in $\mathcal{F}_\Delta$ for $\Delta\geq4$ have maximum degree $\Delta$ and the graphs in $\mathcal{F}_3$ have maximum degree either two or three.

In~\cite{PaperPart1}, we have shown that Conjecture~\ref{conj_G Delta_UB} holds for trees (see Theorem~\ref{thm_trees}). It turns out that, for triangle-free graphs with maximum degree at least $3$, the set of graphs requiring a positive constant $c$ is exactly the same as the set of trees with maximum degree at least $3$ needing a positive constant. Our proof uses the result for trees from~\cite{PaperPart1} as a starting point for an induction. After that, we assume that a triangle-free graph contains at least one cycle containing some edge $e$. We remove that edge to construct a graph $G'$ which, by induction, satisfies the conjecture. Hence, $G'$ contains a \textit{small} identifying code which we can use to construct another identifying code of the same size for $G$. One difficulty for proving the conjecture is the existence of the set of graphs requiring $c>0$. Since $\mathcal F_3$ is the largest among the sets $\mathcal F_\Delta$, the case $\Delta=3$ requires a lot of special argumentation.

\medskip

\textbf{Structure of the paper.} First in Section~\ref{sec:prelim}, we introduce some useful definitions and lemmas. In Subsection~\ref{SecCoiden}, we introduce terminology and results about a generalization of identifying codes, which are later used in the proof of our main result. We continue in Subsection~\ref{secPaths} where we discuss Conjecture~\ref{conj_G Delta_UB} when $\Delta=2$. In Subsection~\ref{secExtremal}, we introduce every connected triangle-free graph which requires a positive constant $c$ for Conjecture~\ref{conj_G Delta_UB}. Understanding these graphs is crucial for the proof of the main theorem. In Section~\ref{sec:main}, we give the proof of Theorem~\ref{thm:trianglefree}. In Section~\ref{sec:beyond}, we use our bound to prove a weaker bound for graphs that have triangles, but can be made triangle-free by removing $t$ edges. Finally, we conclude in Section~\ref{sec:conclu}.

\section{Preliminaries and known results}
\label{sec:prelim}

In the following, we go through the notation we use throughout this article.
We denote by $V(G)$ and $E(G)$ the vertex and edge sets of graph $G=(V(G),E(G))$. We usually denote $n=|V(G)|$. For a set of vertices $S$ we denote $N(S)=\bigcup_{v\in S}N(v)$ and $N[S]=\bigcup_{v\in S}N[v]$. We denote by  $\deg_G(v)=|N_G(v)|$  the \emph{degree} of the vertex $v$ in graph $G$. A \emph{leaf} is a vertex with degree one and its only neighbor  is called \emph{support vertex}. In the literature, a leaf is also known as a \emph{pendant vertex}. Naturally, any vertex of a graph $G$ that is not a leaf of $G$ is  referred to as a \emph{non-leaf} vertex of $G$. We denote the complement of a graph $G$ by $\overline{G}$. We sometimes denote the maximum degree of graph $G$ by $\Delta(G)$, its number of vertices by $n(G)$ and its number of edges by $m(G)$. The \textit{girth} of the graph $G$ refers to the number of vertices in the shortest cycle in graph $G$.

We say that vertex $u$ (or vertex subset $C$) \textit{separates} vertices $v$ and $w$ if vertex $u$ (or some vertex $z \in C$) is in exactly one of sets $N[u]$ and $N[v]$. Given a vertex subset $C$ (often an identifying code), with \textit{code neighborhood} of vertex $u\in V(G)$ we refer to the set $C\cap N[u]$.

On many occasions throughout this article, we look at a subgraph of a graph $G$ obtained by deleting some vertices or edges. To that end, given a graph $G$ and a set $S$ containing some vertices and edges of $G$, we define $G-S$ as the subgraph of $G$ obtained by deleting from $G$ all edges and vertices (and edges incident with the vertices) of $S$.

We use the following lemma multiple times for arguing that some vertices are identified.
\begin{lemma}\label{LemP3}
Let $G$ be a triangle-free graph and $S$, a subset of vertices of $G$. If three vertices $u,v,w$ inducing a $P_3$ are in $S$, then each of them has a unique closed neighborhood within $S$.
\end{lemma}

\begin{proof}
Let  $u,v,w\in S\subseteq V(G)$ induce a $P_3$ in $G$ where $v$ is the middle vertex of the path. Suppose first on the contrary that $N[u]\cap S=N[x]\cap S$ for some $x\in V(G)\setminus\{u\}$. However, now vertices $x,u$ and $v$ form a triangle or if $x=v$, then $x,u$ and $w$ form a triangle. The case with $w$ is symmetric. Suppose then that $N[v]\cap S=N[x]\cap S$ for some $x\in V(G)\setminus\{v\}$. Now, $v,w$ and $x$ form a triangle or if $x=w$, then $u,v$ and $w$ form the triangle. Hence, the claim follows since $G$ is triangle-free.
\end{proof}

\subsection{$(X,Y)$-separating codes and $(X,Y)$-identifying codes}\label{SecCoiden}

We now introduce a generalization of identifying codes, and an upper bound for them that generalizes Bondy's theorem on induced subsets~\cite{B72}, that will be used several times in our proofs and, we believe, can be useful in many other settings as well.

Let $G=(V,E)$ be a graph and $X$ and $Y$ be two (not necessarily disjoint) vertex subsets of $G$. Then, $Y$ induces a partition on $X$ by the equivalence relation $\sim$ defined by $u \sim v$ if and only if $N_G[u] \cap Y = N_G[v] \cap Y$ for any $(u,v) \in X \times X$. If the partition on $X$ induced by $Y$ is such that each part is a singleton set, that is, for each pair $u,v \in X$ there exists a vertex $w \in Y$ that separates $u,v$, then $Y$ is called an \emph{$(X,Y)$-separating set} in $G$, the set $X$ is called $Y$-separable and the graph $G$ is called $(X,Y)$-separable. In addition, any subset $C \subseteq Y$ that is an $(X,C)$-separating set is also referred to as an $(X,Y)$-separating set. Moreover, we call any vertex subset $C$ of $G$  an \emph{$(X,Y)$-identifying code} of $G$ if 1) $C$ is an $(X,Y)$-separating set in $G$, and 2) $C$ is a dominating set of $X$. If such an $(X,Y)$-identifying code of $G$ exists, then  the set $X$ is called \emph{$Y$-identifiable}. Notice that, if $X$ is $Y$-identifiable and $C$ is an $(X,Y)$-identifying code of $G$, then $Y$ itself is an $(X,Y)$-identifying code of $G$. Furthermore, set $X$ is also $C$-identifiable and $C$ is an $(X,C)$-identifying code of $G$. In particular, if $G$ is an identifiable graph and $C$ is an identifying code of $G$, then the vertex set $V$ is $V$-identifiable and $C$-identifiable and $C$ is a $(V,V)$-identifying code and a $(V,C)$-identifying code of $G$. When $X$ and $Y$ are disjoint and induce a bipartite graph, an $(X,Y)$-identifying code has been called a \emph{discriminating code} in the literature~\cite{CCCCHL08}.

\begin{lemma}\label{lemXYID}
Let $G$ be a graph with vertex subsets $X$ and $Y$ such that $X$ is $Y$-identifiable. Then, there is an $(X,Y)$-separating set in $G$ of size at most $|X|-1$, and an $(X,Y)$-identifying code of size at most $|X|$.
\end{lemma}

\begin{proof}
Assume that $G$ is $(X,Y)$-separable. If $|X|=1$, there is nothing to do. Otherwise, we inductively construct an $(X,Y)$-separating code $C$ of $G$ such that $|C|\le |X|-1$. To begin with, let $u,v$ be an arbitrary pair of distinct vertices of $X$ and let $c \in Y$ such that $c$ separates the pair: $c$ exists since $G$ is $(X,Y)$-separable. Then, we let $C = \{c\}$. Let $\mathcal{P}_C$ be the partition induced by $C$ on $X$, where two vertices of $X$ are in the same part if and only if their closed neighborhood in $G$ intersects the same subset of $C$. Then, for as long as there exists a part $P$ of $\mathcal{P}_C$ such that $u',v' \in P$ for two distinct vertices $u',v'$ of $X$, the construction of $C$ follows inductively by selecting an element $c'\in Y$ such that $c'$ separates $u',v'$, and letting $c'\in C$. At each step, since $G$ is $(X,Y)$-separable, such $c'$ exists. Moreover, we notice that at each inductive step, we have $|C| \le |\mathcal{P}_C|-1$, since at each step, we increase the number of parts by at least~1, and the size of $C$ by exactly~1. This implies that we must have $|C| \le |X|-1$, since $|\mathcal{P}_C| \le |X|$, affirming the first claim.

Now, if moreover $G$ is $(X,Y)$-identifiable, we proceed as above to first build the $(X,Y)$-separating set $C$ of $G$ of size at most $|X|-1$. Now, any two vertices of $X$ are separated by $C$. Furthermore, there exists at most one vertex of $X$ that is not dominated by $C$; for otherwise, if there exist two distinct vertices $x,x' \in X$ not dominated by $C$, it would imply that $N_G[x] \cap C = N_G[x'] \cap C = \emptyset$ and so, $C$ would not be an $(X,Y)$-separating set of $X$, a contradiction. Therefore, let $x \in X$ be not dominated by $C$ (if such an $x$ exists). Since $Y$ dominates $X$, there exists $c'' \in N_G[x] \cap Y$. Then, we let $c'' \in C$, thus making $C$ an $(X,Y)$-identifying code of $G$ with $|C|\le |X|$. This completes the proof.
\end{proof}

Lemma~\ref{lemXYID} generalizes Bondy's celebrated theorem on ``induced subsets''~\cite{B72}. Indeed, in Bondy's theorem, one is given a set $X$ of elements and a collection $\mathcal A=\{A_1,\ldots, A_n\}$ of subsets of $X$; it is proved by Bondy that there is a subset of at most $|X|-1$ subsets of $\mathcal A$ that form an $(X,\mathcal A)$-separating set, when viewing $X$ and $\mathcal A$ as the two partite sets of a bipartite graph (the incidence bipartite graph of the hypergraph $(X,\mathcal A)$), provided this graph is $\mathcal A$-identifiable. Bondy's original proof uses an elegant graph-theoretic argument~\cite{B72} (several proofs of algebraic nature have also been provided, see for example~\cite{Winter00}). A similar statement, formulated in the language of graphs, is also proved by Gutin, Ramanujan, Reidl, and Wahlstr\"{o}m in~\cite[Lemma~8]{GRRW2020}, by an inductive argument similar to the one presented here. These prior results however are only concerned with separating sets (thus in their setting, one vertex may remain undominated), and with the special case where $X$ and $Y$ are disjoint. Hence, our result both generalizes and strenghtens these setups.

\subsection{Paths and cycles}\label{secPaths}

Our main result requires a precise understanding of graphs with $\Delta=2$ and triangle-free graphs which need a positive constant $c$ for Conjecture~\ref{conj_G Delta_UB}. Hence, in this subsection, we recall results on all connected graphs with $\Delta=2$, that is, on paths and cycles. A path (cycle) on $n$ vertices is denoted by $P_n$ ($C_n$).

The identification number of all identifiable paths (that is, of all paths except $P_2$) was determined by Bertrand et al.~\cite{BCHL2004}. Moreover, using an upper bound from~\cite{BCHL2004} on even cycles of order at least~$6$, Gravier et al.~\cite{GMS2006} provided the exact values of the identification numbers of all identifiable cycles (that is, cycles of length at least $4$). We summarize  these results in the following theorems.

\begin{theorem}[{\cite[Theorem 3]{BCHL2004}}] \label{thm_BCHL2004}
If $P_n$ is a path on $n$ vertices, then we have
\[
\ID(P_n) =
\begin{cases} \frac{n}{2}+\frac{1}{2}, &\text{if $n \ge 1$ is odd}, \1\\
								  \frac{n}{2}+1, &\text{if $n \ge 4$ is even}.
\end{cases}
\]
\end{theorem}

\begin{theorem}[{\cite[Theorems 2 and 4]{GMS2006}}]\label{thm_GMS2006}
If $C_n$ is a cycle on $n$ vertices, then we have
\[
\ID(C_n) = \begin{cases} 3, &\text{if $n = 4,5$}, \1\\
                                 \frac{n}{2}, &\text{if $n \ge 6$ is even}, \1\\
                                 \frac{n}{2}+\frac{3}{2}, &\text{if $n \ge 7$ is odd}.
\end{cases}
\]
\end{theorem}

Using Theorems~\ref{thm_BCHL2004} and~\ref{thm_GMS2006}, therefore, one has the following corollary.

\begin{corollary} \label{cor_paths & cycles}
The following hold. \\[-16pt]
\begin{enumerate}
\item[{\rm (a)}] If $G$ is a path, then $\ID(P_n) = \lfloor \frac{n}{2} \rfloor+1$. \1
\item[{\rm (b)}] If $G = C_4$ or $G = C_5$, then $\ID(G) = \lfloor \frac{n}{2} \rfloor+1$. \1
\item[{\rm (c)}] If $n \ge 6$ is even, then $\ID(C_n) = \frac{n}{2}$. \1
\item[{\rm (d)}] If $n \ge 7$ is odd, then $\ID(C_n) = \frac{n}{2}+\frac{3}{2}$. \1 
\item[{\rm (e)}] If $n = 4$, then $\ID(P_n) = \frac{3}{4}n$, and if $n \ge 3$ and $n \ne 4$, then $\ID(P_n) \le \frac{2}{3}n$. \1
\item[{\rm (f)}] If $n \in \{4,7\}$, then $\frac{2}{3}n < \ID(C_n) \le \frac{3}{4}n$, and if $n \ge 3$ and $n \notin \{4,7\}$, then $\ID(C_n) \le \frac{2}{3}n$.
\end{enumerate}
\end{corollary}

We shall need the following elementary property of odd cycles with one edge added.

\begin{observation}[Proposition 4.1 of \cite{skaggs2007identifying}]
\label{ob:Ftree1}
If $n=2k+1\ge 7$ is odd, $G = C_n$ and if $e \in E(\overline{G})$, then $\ID(G + e)\le k+1 \le \frac{2}{3}n$.
\end{observation}

\subsection{Extremal triangle-free graphs}\label{secExtremal}

In this section, we define and discuss the exceptional triangle-free graphs of the statement of Theorem~\ref{thm:trianglefree}, that is, those in the set $\mathcal{F}_\Delta$ that require $c>0$ in the bound of Conjecture~\ref{conj_G Delta_UB}. The notation for set $\mathcal{T}_{3}$ (already defined in~\cite{PaperPart1}) will be useful in our proof for the main theorem.

\begin{figure}[!htb]
\centering
\begin{subfigure}[t]{0.3\textwidth}
\centering
\begin{tikzpicture}[
blacknode/.style={circle, draw=black!, fill=black!, thick},
whitenode/.style={circle, draw=black!, fill=white!, thick},
scale=0.5]
\tiny
\node[blacknode] (0) at (0,0) {};
\node[whitenode] (1) at (1.5,0) {};
\node[blacknode] (2) at (2.5,1) {};
\node[blacknode] (3) at (2.5,-1) {};
\draw[-, thick] (0) -- (1);
\draw[-, thick] (1) -- (2);
\draw[-, thick] (1) -- (3);
\end{tikzpicture}
\caption{$T_0$: $\ID(T_0) = 3 = 0.75n$.}
\end{subfigure}
\hspace{2mm}
\begin{subfigure}[t]{0.3\textwidth}
\centering
\begin{tikzpicture}[
blacknode/.style={circle, draw=black!, fill=black!, thick},
whitenode/.style={circle, draw=black!, fill=white!, thick},
scale=0.5]
\tiny
\node[blacknode] (0) at (0,0) {};
\node[whitenode] (1) at (1.5,0) {};
\node[blacknode] (2) at (2.5,1) {};
\node[blacknode] (3) at (2.5,-1) {};
\node[whitenode] (1') at (-1.5,0) {};
\node[blacknode] (2') at (-2.5,1) {};
\node[blacknode] (3') at (-2.5,-1) {};
\draw[-, thick] (0) -- (1);
\draw[-, thick] (1) -- (2);
\draw[-, thick] (1) -- (3);
\draw[-, thick] (0) -- (1');
\draw[-, thick] (1') -- (2');
\draw[-, thick] (1') -- (3');
\end{tikzpicture}
\caption{$T_1$: $\ID(T_1) = 5 \approx 0.71n$.}
\end{subfigure}
\hspace{2mm}
\begin{subfigure}[t]{0.3\textwidth}
\centering
\begin{tikzpicture}[
blacknode/.style={circle, draw=black!, fill=black!, thick},
whitenode/.style={circle, draw=black!, fill=white!, thick},
scale=0.5]
\tiny
\node[whitenode] (1) at (0,0) {};
\node[blacknode] (2) at (-1,1) {};
\node[blacknode] (3) at (-1,-1) {};
\node[blacknode] (4) at (1.5,0) {};
\node[blacknode] (5) at (3,0) {};
\node[blacknode] (6) at (4.5,0) {};
\node[whitenode] (7) at (6,0) {};
\draw[-, thick, black!] (1) -- (2);
\draw[-, thick, black!] (1) -- (3);
\draw[-, thick, black!] (1) -- (4);
\draw[-, thick] (5) -- (4);
\draw[-, thick] (4) -- (6);
\draw[-, thick] (6) -- (7);
\end{tikzpicture}
\caption{$T_2$: $\ID(T_2) = 5 \approx 0.71n$.}
\end{subfigure}\vspace{8mm}
\hspace{2mm}
\begin{subfigure}[t]{0.3\textwidth}
\centering
\begin{tikzpicture}[
blacknode/.style={circle, draw=black!, fill=black!, thick},
whitenode/.style={circle, draw=black!, fill=white!, thick},
scale=0.5]
\tiny
\node[whitenode] (1) at (-1,1) {};
\node[blacknode] (2) at (-1,2.5) {};
\node[blacknode] (3) at (-2.5,1) {};
\node[whitenode] (4) at (0,0) {};
\node[whitenode] (8) at (-1,-1) {};
\node[blacknode] (9) at (-1,-2.5) {};
\node[blacknode] (10) at (-2.5,-1) {};
\node[blacknode] (5) at (1.5,0) {};
\node[blacknode] (6) at (3,0) {};
\node[blacknode] (7) at (4.5,0) {};
\draw[-, thick, black!] (1) -- (2);
\draw[-, thick, black!] (1) -- (3);
\draw[-, thick, black!] (1) -- (4);
\draw[-, thick, black!] (4) -- (8);
\draw[-, thick, black!] (8) -- (9);
\draw[-, thick, black!] (8) -- (10);
\draw[-, thick] (5) -- (4);
\draw[-, thick] (4) -- (6);
\draw[-, thick] (6) -- (7);
\end{tikzpicture}
\caption{$T_3$: $\ID(T_3) = 7 = 0.7n$.}
\end{subfigure}
\hspace{2mm}
\begin{subfigure}[t]{0.3\textwidth}
\centering
\begin{tikzpicture}[
blacknode/.style={circle, draw=black!, fill=black!, thick},
whitenode/.style={circle, draw=black!, fill=white!, thick},
scale=0.5]
\tiny
\node[blacknode] (0) at (0,0) {};
\node[whitenode] (1) at (1.5,0) {};
\node[blacknode] (2) at (2.5,1) {};
\node[blacknode] (3) at (2.5,-1) {};
\node[whitenode] (4) at (-1,-1) {};
\node[blacknode] (5) at (-1,-2.5) {};
\node[blacknode] (6) at (-2.6,-1) {};
\node[whitenode] (7) at (-1,1) {};
\node[blacknode] (8) at (-1,2.5) {};
\node[blacknode] (9) at (-2.6,1) {};
\draw[-, thick] (0) -- (1);
\draw[-, thick] (1) -- (2);
\draw[-, thick] (1) -- (3);
\draw[-, thick] (0) -- (4);
\draw[-, thick] (4) -- (5);
\draw[-, thick] (4) -- (6);
\draw[-, thick] (0) -- (7);
\draw[-, thick] (7) -- (8);
\draw[-, thick] (7) -- (9);
\end{tikzpicture}
\caption{$T_4$: $\ID(T_4) = 7 = 0.7n$.}
\end{subfigure}
\hspace{0mm}
\begin{subfigure}[t]{0.3\textwidth}
\centering
\begin{tikzpicture}[
blacknode/.style={circle, draw=black!, fill=black!, thick},
whitenode/.style={circle, draw=black!, fill=white!, thick},
scale=0.5]
\tiny
\node[blacknode] (0) at (0,0) {};
\node[whitenode] (1) at (1.5,0) {};
\node[blacknode] (2) at (2.5,1) {};
\node[blacknode] (3) at (2.5,-1) {};
\node[whitenode] (13) at (4,1.2) {};
\node[blacknode] (14) at (5.2,2.4) {};
\node[blacknode] (15) at (5.4,0.5) {};
\node[whitenode] (1') at (-1.5,0) {};
\node[blacknode] (2') at (-2.5,1) {};
\node[blacknode] (3') at (-2.5,-1) {};
\draw[-, thick] (0) -- (1);
\draw[-, thick] (1) -- (2);
\draw[-, thick] (1) -- (3);
\draw[-, thick] (2) -- (13);
\draw[-, thick] (13) -- (14);
\draw[-, thick] (13) -- (15);
\draw[-, thick] (0) -- (1');
\draw[-, thick] (1') -- (2');
\draw[-, thick] (1') -- (3');
\end{tikzpicture}
\caption{$T_5$: $\ID(T_5) = 7 = 0.7n$.}
\end{subfigure}\vspace{8mm}
\hspace{2mm}
\begin{subfigure}[t]{0.3\textwidth}
\centering
\begin{tikzpicture}[
blacknode/.style={circle, draw=black!, fill=black!, thick},
whitenode/.style={circle, draw=black!, fill=white!, thick},
scale=0.5]
\tiny
\node[blacknode] (0) at (0,0) {};
\node[whitenode] (1) at (1.5,0) {};
\node[blacknode] (2) at (2.5,1) {};
\node[blacknode] (3) at (2.5,-1) {};
\node[whitenode] (13) at (4,1.2) {};
\node[blacknode] (14) at (5.2,2.4) {};
\node[blacknode] (15) at (5.4,0.5) {};
\node[whitenode] (19) at (4,-1.2) {};
\node[blacknode] (20) at (5.2,-2.4) {};
\node[blacknode] (21) at (5.4,-0.5) {};
\node[whitenode] (1') at (-1.5,0) {};
\node[blacknode] (2') at (-2.5,1) {};
\node[blacknode] (3') at (-2.5,-1) {};
\draw[-, thick] (0) -- (1);
\draw[-, thick] (1) -- (2);
\draw[-, thick] (1) -- (3);
\draw[-, thick] (2) -- (13);
\draw[-, thick] (13) -- (14);
\draw[-, thick] (13) -- (15);
\draw[-, thick] (3) -- (19);
\draw[-, thick] (19) -- (20);
\draw[-, thick] (19) -- (21);
\draw[-, thick] (0) -- (1');
\draw[-, thick] (1') -- (2');
\draw[-, thick] (1') -- (3');
\end{tikzpicture}
\caption{$T_6$: $\ID(T_6) = 9 \approx 0.69n$.}
\end{subfigure}
\hspace{2mm}
\begin{subfigure}[t]{0.3\textwidth}
\centering
\begin{tikzpicture}[
blacknode/.style={circle, draw=black!, fill=black!, thick},
whitenode/.style={circle, draw=black!, fill=white!, thick},
scale=0.5]
\tiny
\node[blacknode] (0) at (0,0) {};
\node[whitenode] (1) at (1.5,0) {};
\node[blacknode] (2) at (2.5,1) {};
\node[blacknode] (3) at (2.5,-1) {};
\node[whitenode] (4) at (-1,-1) {};
\node[blacknode] (5) at (-1,-2.5) {};
\node[blacknode] (6) at (-2.5,-1) {};
\node[whitenode] (7) at (-1,1) {};
\node[blacknode] (8) at (-1,2.5) {};
\node[blacknode] (9) at (-2.5,1) {};
\node[whitenode] (10) at (2.5,2.5) {};
\node[blacknode] (11) at (3.5,3.5) {};
\node[blacknode] (12) at (1.5,3.5) {};
\draw[-, thick] (0) -- (1);
\draw[-, thick] (1) -- (2);
\draw[-, thick] (1) -- (3);
\draw[-, thick] (2) -- (10);
\draw[-, thick] (10) -- (11);
\draw[-, thick] (10) -- (12);
\draw[-, thick] (0) -- (4);
\draw[-, thick] (4) -- (5);
\draw[-, thick] (4) -- (6);
\draw[-, thick] (0) -- (7);
\draw[-, thick] (7) -- (8);
\draw[-, thick] (7) -- (9);
\end{tikzpicture}
\caption{$T_7$: $\ID(T_7) = 9 \approx 0.69n$.}
\end{subfigure}
\hspace{0mm}
\begin{subfigure}[t]{0.3\textwidth}
\centering
\begin{tikzpicture}[
blacknode/.style={circle, draw=black!, fill=black!, thick},
whitenode/.style={circle, draw=black!, fill=white!, thick},
scale=0.5]
\tiny
\node[blacknode] (0) at (0,0) {};
\node[whitenode] (1) at (1.5,0) {};
\node[blacknode] (2) at (2.5,1) {};
\node[blacknode] (3) at (2.5,-1) {};
\node[whitenode] (13) at (4,1.2) {};
\node[blacknode] (14) at (5.2,2.4) {};
\node[blacknode] (15) at (5.4,0.5) {};
\node[whitenode] (19) at (4,-1.2) {};
\node[blacknode] (20) at (5.2,-2.4) {};
\node[blacknode] (21) at (5.4,-0.5) {};
\node[whitenode] (4) at (-1,-1) {};
\node[blacknode] (5) at (-1,-2.5) {};
\node[blacknode] (6) at (-2.5,-1) {};
\node[whitenode] (7) at (-1,1) {};
\node[blacknode] (8) at (-1,2.5) {};
\node[blacknode] (9) at (-2.5,1) {};
\draw[-, thick] (0) -- (1);
\draw[-, thick] (1) -- (2);
\draw[-, thick] (1) -- (3);
\draw[-, thick] (2) -- (13);
\draw[-, thick] (13) -- (14);
\draw[-, thick] (13) -- (15);
\draw[-, thick] (3) -- (19);
\draw[-, thick] (19) -- (20);
\draw[-, thick] (19) -- (21);
\draw[-, thick] (0) -- (4);
\draw[-, thick] (4) -- (5);
\draw[-, thick] (4) -- (6);
\draw[-, thick] (0) -- (7);
\draw[-, thick] (7) -- (8);
\draw[-, thick] (7) -- (9);
\end{tikzpicture}
\caption{$T_8$: $\ID(T_8) = 11 \approx 0.69n$.}
\end{subfigure}\vspace{8mm}
\hspace{6mm}
\begin{subfigure}[t]{0.3\textwidth}
\centering
\begin{tikzpicture}[
blacknode/.style={circle, draw=black!, fill=black!, thick},
whitenode/.style={circle, draw=black!, fill=white!, thick},
scale=0.5]
\tiny
\node[blacknode] (0) at (0,0) {};
\node[whitenode] (1) at (1.5,0) {};
\node[blacknode] (2) at (2.5,1) {};
\node[blacknode] (3) at (2.5,-1) {};
\node[whitenode] (4) at (-1,-1) {};
\node[blacknode] (5) at (-1,-2.5) {};
\node[blacknode] (6) at (-2.5,-1) {};
\node[whitenode] (7) at (-1,1) {};
\node[blacknode] (8) at (-1,2.5) {};
\node[blacknode] (9) at (-2.5,1) {};
\node[whitenode] (10) at (2.5,2.5) {};
\node[blacknode] (11) at (3.5,3.5) {};
\node[blacknode] (12) at (1.5,3.5) {};
\node[whitenode] (13) at (4,1.2) {};
\node[blacknode] (14) at (5.2,2.4) {};
\node[blacknode] (15) at (5.4,0.5) {};
\draw[-, thick] (0) -- (1);
\draw[-, thick] (1) -- (2);
\draw[-, thick] (1) -- (3);
\draw[-, thick] (2) -- (10);
\draw[-, thick] (10) -- (11);
\draw[-, thick] (10) -- (12);
\draw[-, thick] (2) -- (13);
\draw[-, thick] (13) -- (14);
\draw[-, thick] (13) -- (15);
\draw[-, thick] (0) -- (4);
\draw[-, thick] (4) -- (5);
\draw[-, thick] (4) -- (6);
\draw[-, thick] (0) -- (7);
\draw[-, thick] (7) -- (8);
\draw[-, thick] (7) -- (9);

\end{tikzpicture}
\caption{$T_9$: $\ID(T_9) = 11 \approx 0.69n$.}
\end{subfigure}
\hspace{4mm}
\begin{subfigure}[t]{0.3\textwidth}
\centering
\begin{tikzpicture}[
blacknode/.style={circle, draw=black!, fill=black!, thick},
whitenode/.style={circle, draw=black!, fill=white!, thick},
scale=0.5]
\tiny
\node[blacknode] (0) at (0,0) {};
\node[whitenode] (1) at (1.5,0) {};
\node[blacknode] (2) at (2.5,1) {};
\node[blacknode] (3) at (2.5,-1) {};
\node[whitenode] (4) at (-1,-1) {};
\node[blacknode] (5) at (-1,-2.5) {};
\node[blacknode] (6) at (-2.5,-1) {};
\node[whitenode] (7) at (-1,1) {};
\node[blacknode] (8) at (-1,2.5) {};
\node[blacknode] (9) at (-2.5,1) {};
\node[whitenode] (10) at (2.5,2.5) {};
\node[blacknode] (11) at (3.5,3.5) {};
\node[blacknode] (12) at (1.5,3.5) {};
\node[whitenode] (13) at (4,1.2) {};
\node[blacknode] (14) at (5.2,2.4) {};
\node[blacknode] (15) at (5.4,0.5) {};
\node[whitenode] (16) at (2.5,-2.5) {};
\node[blacknode] (17) at (3.5,-3.5) {};
\node[blacknode] (18) at (1.5,-3.5) {};
\draw[-, thick] (0) -- (1);
\draw[-, thick] (1) -- (2);
\draw[-, thick] (1) -- (3);
\draw[-, thick] (2) -- (10);
\draw[-, thick] (10) -- (11);
\draw[-, thick] (10) -- (12);
\draw[-, thick] (2) -- (13);
\draw[-, thick] (13) -- (14);
\draw[-, thick] (13) -- (15);
\draw[-, thick] (3) -- (16);
\draw[-, thick] (16) -- (17);
\draw[-, thick] (16) -- (18);
\draw[-, thick] (0) -- (4);
\draw[-, thick] (4) -- (5);
\draw[-, thick] (4) -- (6);
\draw[-, thick] (0) -- (7);
\draw[-, thick] (7) -- (8);
\draw[-, thick] (7) -- (9);
\end{tikzpicture}
\caption{$T_{10}$: $\ID(T_{10}) = 13 \approx 0.68n$.}
\end{subfigure}
\hspace{4mm}
\begin{subfigure}[t]{0.3\textwidth}
\centering
\begin{tikzpicture}[
blacknode/.style={circle, draw=black!, fill=black!, thick},
whitenode/.style={circle, draw=black!, fill=white!, thick},
scale=0.5]
\tiny
\node[blacknode] (0) at (0,0) {};
\node[whitenode] (1) at (1.5,0) {};
\node[blacknode] (2) at (2.5,1) {};
\node[blacknode] (3) at (2.5,-1) {};
\node[whitenode] (4) at (-1,-1) {};
\node[blacknode] (5) at (-1,-2.5) {};
\node[blacknode] (6) at (-2.5,-1) {};
\node[whitenode] (7) at (-1,1) {};
\node[blacknode] (8) at (-1,2.5) {};
\node[blacknode] (9) at (-2.5,1) {};
\node[whitenode] (10) at (2.5,2.5) {};
\node[blacknode] (11) at (3.5,3.5) {};
\node[blacknode] (12) at (1.5,3.5) {};
\node[whitenode] (13) at (4,1.2) {};
\node[blacknode] (14) at (5.2,2.4) {};
\node[blacknode] (15) at (5.4,0.5) {};
\node[whitenode] (16) at (2.5,-2.5) {};
\node[blacknode] (17) at (3.5,-3.5) {};
\node[blacknode] (18) at (1.5,-3.5) {};
\node[whitenode] (19) at (4,-1.2) {};
\node[blacknode] (20) at (5.2,-2.4) {};
\node[blacknode] (21) at (5.4,-0.5) {};
\draw[-, thick] (0) -- (1);
\draw[-, thick] (1) -- (2);
\draw[-, thick] (1) -- (3);
\draw[-, thick] (2) -- (10);
\draw[-, thick] (10) -- (11);
\draw[-, thick] (10) -- (12);
\draw[-, thick] (2) -- (13);
\draw[-, thick] (13) -- (14);
\draw[-, thick] (13) -- (15);
\draw[-, thick] (3) -- (16);
\draw[-, thick] (16) -- (17);
\draw[-, thick] (16) -- (18);
\draw[-, thick] (3) -- (19);
\draw[-, thick] (19) -- (20);
\draw[-, thick] (19) -- (21);
\draw[-, thick] (0) -- (4);
\draw[-, thick] (4) -- (5);
\draw[-, thick] (4) -- (6);
\draw[-, thick] (0) -- (7);
\draw[-, thick] (7) -- (8);
\draw[-, thick] (7) -- (9);
\end{tikzpicture}
\caption{$T_{11}$: $\ID(T_{11}) = 15 \approx 0.68n$.}
\end{subfigure}
\caption{The family~$\mathcal{T}_{3}$ of trees of maximum degree~3 requiring $c>0$ in Conjecture~\ref{conj_G Delta_UB}. The set of black vertices in each figure constitutes an identifying code of the tree.}
\label{fig:trees}
\end{figure}
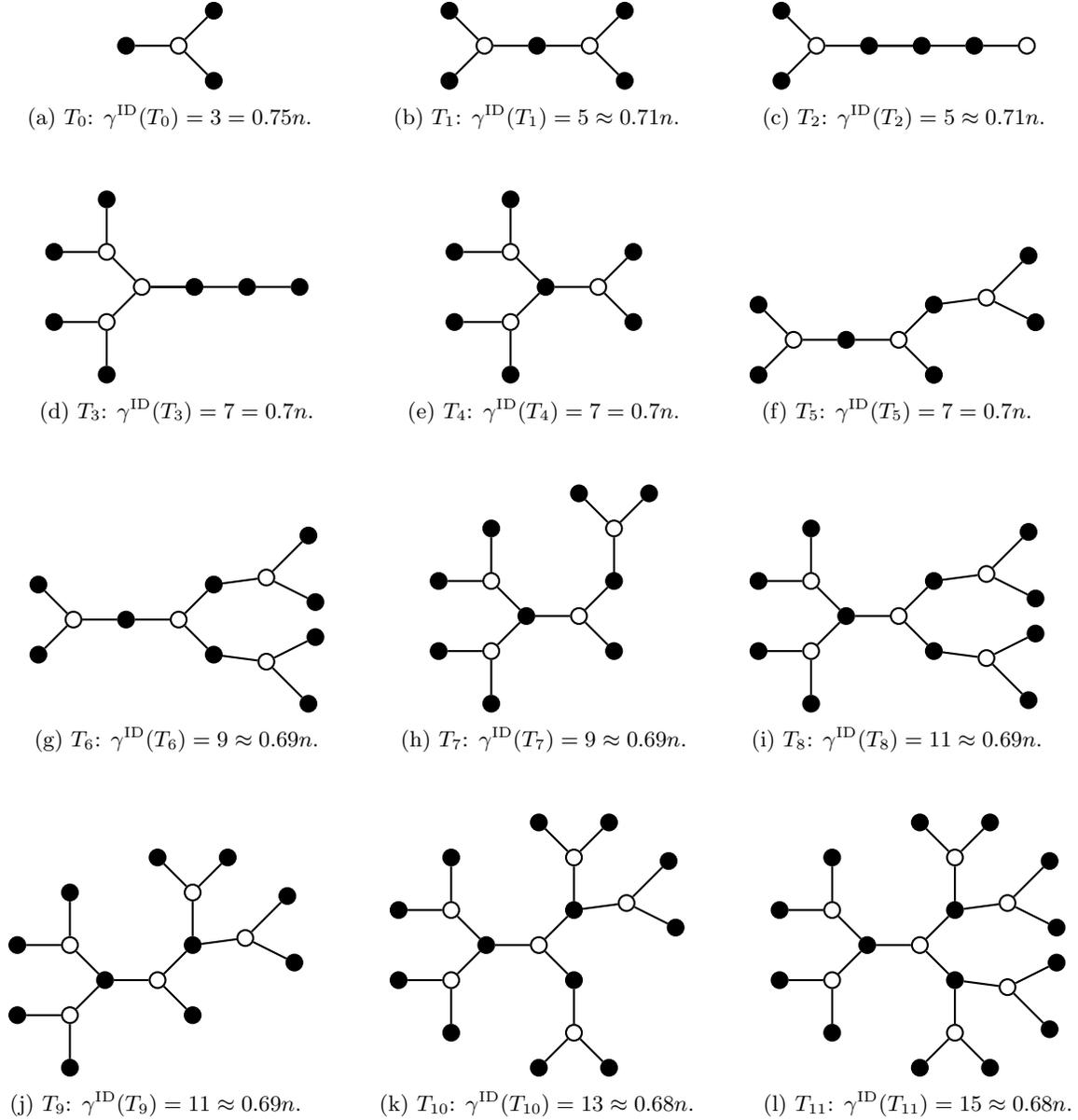

\begin{definition}\label{def:TreeFDelta}
For $\Delta = 3$, we define $\mathcal{T}_{3} = \{T_0,T_1,T_2,\ldots,T_{11}\}$ to be the collection of 12 trees of maximum degree~$3$ as in Figure~\ref{fig:trees}, and for $\Delta \ge 4$, we let $\mathcal{T}_\Delta = \{K_{1,\Delta} \}$.

For $\Delta = 3$, we let $\mathcal{F}_{3} = \mathcal{T}_{3} \cup \{P_4, C_4, C_7\}$ and for $\Delta \ge 4$, we let $\mathcal{F}_{\Delta} = \mathcal{T}_\Delta = \{K_{1,\Delta}\}$.
\end{definition}

We note that $\ID(T_0) = 3$, $\ID(T_1) = \ID(T_2) = 5$,  $\ID(T_3) = \ID(T_4) = \ID(T_5) = 7$, $\ID(T_6) = \ID(T_7) = 9$, $\ID(T_8) = \ID(T_9) = 11$, $\ID(T_{10}) = 13$, and $\ID(T_{11}) = 15$. Generally we have the following.

\begin{proposition}[\cite{PaperPart1}]\label{prop-T_Delta}
If $\Delta\ge 3$ is an integer and $T$ is a tree of order $n$ in $\mathcal T_\Delta$, then $\ID(T)= \left( \frac{\Delta - 1}{\Delta}\right) n+\frac{1}{\Delta}$.
\end{proposition}

By Corollary~\ref{cor_paths & cycles}, if $G$ has maximum degree $\Delta=2$, then $\ID(G)\le \left( \frac{\Delta - 1}{\Delta}\right) n+\frac{3}{2}$. When $\Delta\ge 3$ we have the following.

\begin{proposition}\label{prop-F_Delta}
If $\Delta\ge 3$ is an integer and $G$ is a graph of order $n$ and maximum degree at most $\Delta$ in $\mathcal F_\Delta$ (possibly, if $\Delta=3$, the maximum degree of $G$ is~2), then $\ID(G)= \left( \frac{\Delta - 1}{\Delta}\right) n+\frac{1}{\Delta}$.
\end{proposition}
\begin{proof}
If $G$ is a tree in $\mathcal T_\Delta$, then this follows from Proposition~\ref{prop-T_Delta}. Otherwise, $\Delta=3$ and $G\in\{P_4,C_4,C_7\}$. If $G\in\{P_4,C_4\}$, by Corollary~\ref{cor_paths & cycles}, $\ID(G)=3=\frac{2}{3}n+\frac{1}{3}$ and if $G=C_7$, $\ID(G)=5=\frac{2}{3}n+\frac{1}{3}$.
\end{proof}

We have shown that trees satisfy Conjecture~\ref{conj_G Delta_UB} in~\cite{PaperPart1}, as follows.
\begin{theorem}[\cite{PaperPart1}]
\label{thm_trees}
If $T\notin \mathcal{T}_{\Delta}$ is a tree of order $n \ge 3$ with maximum degree $\Delta \ge 3$, then
\[
\ID(T) \le\left( \frac{\Delta - 1}{\Delta}\right) n.
\]

\end{theorem}

 In the following proposition, we present some useful properties for the structure of identifying codes in trees in set $\mathcal{T}_{3}$. The proof is provided in
 ~\cite{PaperPart1}.

\begin{proposition}[\cite{PaperPart1}]\label{propTreeCodes}
If $T$ is a tree in $\mathcal{T}_3$, then the following properties hold. \\ [-20pt]
\begin{enumerate}
    \item[(i)] If $T\neq T_2$, then $T$ has an optimal identifying code $C(T)$ which includes all vertices of degree at most~2.
    \item[(ii)] If $T\notin\{T_2,T_3\}$, then $C(T)$ can be chosen as an independent set. When we delete any code vertex $v$ from $T$, set $C(T)\setminus\{v\}$ forms an optimal identifying code of the forest $T-v$.
\end{enumerate}
\end{proposition}

We shall also need the following property of trees in the family $\mathcal{T}_{3}$.

\begin{lemma}
\label{ob:treeFam}
Let $T \in \mathcal{T}_{3}$ be a tree of order~$n$.
If $e \in E(\overline{T})$ is such that $T + e$ is triangle-free and $\Delta(T + e) = 3$, then $\ID(T + e) < \frac{2}{3}n$.
\end{lemma}

\begin{proof}
Assume first that $T\not\in\{T_2,T_3\}$. Note that since $T+e$ is triangle-free, $T\neq T_0$.
Let $C(T)$ be an optimal identifying code in $T$ which includes every vertex of degree at most $2$ such that $C(T)$ is also an independent set. Such an identifying code exists by Proposition~\ref{propTreeCodes} and is presented in Figure~\ref{fig:trees}. Observe that $C(T)$ 3-dominates every vertex in $V(T)\setminus C(T)$. Let us denote the endpoints of edge $e$ by $u$ and $v$. Since $\Delta(T+e)=3$, both $u$ and $v$ have degree at most 2 in $T$ and hence, $u,v\in C(T)$. Furthermore, since $T+e$ is triangle-free, at least one of $u,v$ is a leaf (since in all trees in $\mathcal T_3\setminus\{T_2,T_3\}$, all degree~2 vertices are at distance~2 apart). Let us assume without loss of generality that $u$ is a leaf and $u_s$ is the adjacent support vertex. Furthermore, let $u'\in C(T)$ be the vertex at distance~2 from $u$ with $d_T(u,v)=d_T(u',v)$ (vertex $u'$ exists due to the structure of trees in $\mathcal T_3\setminus\{T_2,T_3\}$ and the assumption that $C(T)$ 3-dominates every vertex in $V(T)\setminus C(T)$).

Consider the set $C=(C(T)\setminus\{v,u'\})\cup\{u_s\}$. Since $|C|=|C(T)|-1$, we have $|C|< \frac{2}{3}n$. We claim that $C$ is an identifying code in $T+e$. Since each vertex in $V(T)\setminus C(T)$ is 3-dominated by $C(T)$, each vertex in $V(T)\setminus C(T)$ is still dominated by $C$. Furthermore, $u'$ is dominated by $u_s$ and $v$ is dominated by $u$. Thus, $C$ is a dominating set in $T+e$. We now check that $C$ is separating. Since only $u'$ and $v$ have been removed from $C(T)$ when constructing $C$, it suffices to show that all vertices in $N[u']\cup N[v]$ are separated from all other vertices of $T+e$, as well as $u$ and $v$ are pairwise separated since we added the edge $e$. (Other pairs of vertices, such as the neighbors of $u_s$, may have been affected, but since we only added code vertices to their closed neighborhood, they remain separated by $C$.) We note that there exists a vertex $u_c\in C\cap N(u_s)\setminus \{u\}$ since $u_s$ is 3-dominated by $C(T)$. Hence, $u,u_s$ and $u_c$ form a 3-path in $C$ and they have a unique neighborhood in $C$ by Lemma~\ref{LemP3}. Furthermore, vertex $v$ is the only vertex with exactly $u$ in its code neighborhood while any neighbor of $v$ is in $V(T)\setminus C$ and remains at 2-dominated or 3-dominated. Hence, vertices in $N(v)$ have unique code neighborhoods. Vertex $u'$ is the only vertex with exactly $u_s$ in its code neighborhood and any vertex in $N(u')\setminus C$ has a unique code neighborhood since it is at least 2-dominated. Therefore, $C$ is an identifying code in $T+e$ satisfying the claim. Hence, we are only left with considering trees $T_2$ and $T_3$.

For the case $T=T_2$, we have presented each possible triangle-free graph $T+e$ of maximum degree~3 in Figures~\ref{fig:T23e}(a)-(d), together with identifying codes satisfying the claim. For $T=T_3$, notice that we may obtain $T_3$ from $T_2$ by attaching a $P_3$ from the middle vertex to a vertex of $T_2$. Hence, when the new edge $e$ is not adjacent to one of the vertices in the $P_3$, we may use the identifying code in $T_2+e$ together with the two new leaves in the newly attached $P_3$. When we have the edge $e$ adjacent to the vertices in the newly attached $P_3$, we may use the identifying code from Figure~\ref{fig:T23e}(e).
\end{proof}

\begin{figure}
\centering
\begin{subfigure}[t]{0.3\textwidth}
\centering
\begin{tikzpicture}[
blacknode/.style={circle, draw=black!, fill=black!, thick},
whitenode/.style={circle, draw=black!, fill=white!, thick},
scale=0.5]
\tiny

\node[blacknode] (1) at (0,0) {};
\node[blacknode] (2) at (-1,1) {};
\node[whitenode] (3) at (-1,-1) {};
\node[blacknode] (4) at (1.5,0) {};
\node[blacknode] (5) at (3,0) {};
\node[whitenode] (6) at (4.5,0) {};
\node[whitenode] (7) at (6,0) {};
\draw[-, thick, black!] (1) -- (2);
\draw[-, thick, black!] (1) -- (3);
\draw[-, thick, black!] (1) -- (4);
\draw[-, thick] (5) -- (4);
\draw[-, thick] (4) -- (6);
\draw[-, thick] (6) -- (7);
\draw[-, thick] (2) to[bend left=10,looseness=1] (7);
\end{tikzpicture}
\caption{$T_2+e_1$.}
\end{subfigure}\qquad
\begin{subfigure}[t]{0.3\textwidth}
\centering
\begin{tikzpicture}[
blacknode/.style={circle, draw=black!, fill=black!, thick},
whitenode/.style={circle, draw=black!, fill=white!, thick},
scale=0.5]
\tiny

\node[blacknode] (1) at (0,0) {};
\node[whitenode] (2) at (-1,1) {};
\node[whitenode] (3) at (-1,-1) {};
\node[blacknode] (4) at (1.5,0) {};
\node[blacknode] (5) at (3,0) {};
\node[blacknode] (6) at (4.5,0) {};
\node[whitenode] (7) at (6,0) {};
\draw[-, thick, black!] (1) -- (2);
\draw[-, thick, black!] (1) -- (3);
\draw[-, thick, black!] (1) -- (4);
\draw[-, thick] (5) -- (4);
\draw[-, thick] (4) -- (6);
\draw[-, thick] (6) -- (7);
\draw[-, thick] (2) to[bend left=10,looseness=1] (6);
\end{tikzpicture}
\caption{$T_2+e_2$.}
\end{subfigure}\qquad
\begin{subfigure}[t]{0.3\textwidth}
\centering
\begin{tikzpicture}[
blacknode/.style={circle, draw=black!, fill=black!, thick},
whitenode/.style={circle, draw=black!, fill=white!, thick},
scale=0.5]
\tiny

\node[blacknode] (1) at (0,0) {};
\node[whitenode] (2) at (-1,1) {};
\node[whitenode] (3) at (-1,-1) {};
\node[blacknode] (4) at (1.5,0) {};
\node[blacknode] (5) at (3,0) {};
\node[blacknode] (6) at (4.5,0) {};
\node[whitenode] (7) at (6,0) {};
\draw[-, thick, black!] (1) -- (2);
\draw[-, thick, black!] (1) -- (3);
\draw[-, thick, black!] (1) -- (4);
\draw[-, thick] (5) -- (4);
\draw[-, thick] (4) -- (6);
\draw[-, thick] (6) -- (7);
\draw[-, thick] (2) to[bend left=10,looseness=1] (5);
\end{tikzpicture}
\caption{$T_2+e_3$.}
\end{subfigure}\\[5mm]
\begin{subfigure}[t]{0.3\textwidth}
\centering
\begin{tikzpicture}[
blacknode/.style={circle, draw=black!, fill=black!, thick},
whitenode/.style={circle, draw=black!, fill=white!, thick},
scale=0.5]
\tiny

\node[blacknode] (1) at (0,0) {};
\node[blacknode] (2) at (-1,1) {};
\node[whitenode] (3) at (-1,-1) {};
\node[blacknode] (4) at (1.5,0) {};
\node[blacknode] (5) at (3,0) {};
\node[whitenode] (6) at (4.5,0) {};
\node[whitenode] (7) at (6,0) {};
\draw[-, thick, black!] (1) -- (2);
\draw[-, thick, black!] (1) -- (3);
\draw[-, thick, black!] (1) -- (4);
\draw[-, thick] (5) -- (4);
\draw[-, thick] (4) -- (6);
\draw[-, thick] (6) -- (7);
\draw[-, thick] (4) to[bend left=30,looseness=1] (7);
\end{tikzpicture}
\caption{$T_2+e_4$.}
\end{subfigure}\qquad
\begin{subfigure}[t]{0.3\textwidth}
\centering
\begin{tikzpicture}[
blacknode/.style={circle, draw=black!, fill=black!, thick},
whitenode/.style={circle, draw=black!, fill=white!, thick},
scale=0.5]
\tiny
\node[blacknode] (1) at (-1,1) {};
\node[whitenode] (2) at (-1,2.5) {};
\node[blacknode] (3) at (-2.5,1) {};
\node[blacknode] (4) at (0,0) {};
\node[whitenode] (8) at (-1,-1) {};
\node[blacknode] (9) at (-1,-2.5) {};
\node[whitenode] (10) at (-2.5,-1) {};
\node[blacknode] (5) at (1.5,0) {};
\node[blacknode] (6) at (3,0) {};
\node[whitenode] (7) at (4.5,0) {};
\draw[-, thick, black!] (1) -- (2);
\draw[-, thick, black!] (1) -- (3);
\draw[-, thick, black!] (1) -- (4);
\draw[-, thick, black!] (4) -- (8);
\draw[-, thick, black!] (8) -- (9);
\draw[-, thick, black!] (8) -- (10);
\draw[-, thick] (5) -- (4);
\draw[-, thick] (4) -- (6);
\draw[-, thick] (6) -- (7);
\draw[-, thick] (3) -- (10);
\end{tikzpicture}
\caption{$T_3+e$.}
\end{subfigure}
    \caption{Some triangle-free graphs of types $T_2+e$ and $T_3+e$.}
    \label{fig:T23e}
\end{figure}

We have listed in Table~\ref{TableConstVal} the graphs which are known (to us) to require a positive constant for Conjecture~\ref{conj_G Delta_UB}. Among these, are the extremal graphs discovered in~\cite{FGKNPV11} (those graphs of order $n$ with identifying code number $n-1$). Those are either stars, or can be built from any number of complements of half-graphs\footnote{A \emph{half-graph} is a special bipartite graph with both parts of the same size, where each part can be ordered so that the open neighborhoods of consecutive vertices differ by exactly one vertex~\cite{EH84}. Their complements can also be described as powers of paths of the form $P_{2k}^{k-1}$~\cite{FGKNPV11}.} by taking their complete join, and optionally, adding a single universal vertex. Note that the latter examples have large cliques and thus, are far from triangle-free. They have maximum degree $n-1$ or $n-2$, and so, they need $c=1/\Delta$ or $c=2/\Delta$ in the bound of Conjecture~\ref{conj_G Delta_UB}.

It is an interesting open question whether there exist any other such graphs and if the constant $c=\frac{3}{2}$ is enough for all graphs. Notice that by Theorem~\ref{thm:trianglefree}, if there exists a graph not listed in Table~\ref{TableConstVal} that requires a positive constant $c$, then it must contain triangles.
\begin{table}[h]
\centering
\begin{tabular}{|c|c|c|c|}
\hline
Graph class & $\Delta$  & $c$  & Reference\\ \hline
\hline
Odd paths & 2  & $1/\Delta=1/2$  & \cite{BCHL2004} (Theorem~\ref{thm_BCHL2004}) \\ \hline
Even paths & 2 & $ 2/\Delta=1$~~~  & \cite{BCHL2004} (Theorem~\ref{thm_BCHL2004})\\ \hline
Odd cycles $C_n$ for $n\geq7$ & 2  & $ 3/\Delta=3/2$  & \cite{GMS2006} (Theorem~\ref{thm_GMS2006})\\ \hline
$C_4$ & 2  & $ 2/\Delta=1$~~~  & \cite{GMS2006} (Theorem~\ref{thm_GMS2006})\\ \hline
$C_5$ & 2  & $ 1/\Delta=1/2$  & \cite{GMS2006} (Theorem~\ref{thm_GMS2006})\\
\hline
$\mathcal{T}_{3}$ & 3  & $ 1/\Delta=1/3$  & \cite{PaperPart1} \\ \hline
$K_{1,\Delta}$ & $\Delta \ge 3$  & $ 1/\Delta$  &  \cite{PaperPart1}\\ \hline
\makecell{Complements of half-graphs \\ and their complete joins} & even $\Delta\ge 2$ & $ 2/\Delta$  & \cite{FGKNPV11} \\ \hline
\makecell{Complements of half-graphs \\ and their complete joins \\ plus one universal vertex} &  odd $\Delta\ge 3$  & $ 1/\Delta$  & \cite{FGKNPV11} \\ \hline
\end{tabular}\caption{Known graphs requiring a positive constant $c$ for Conjecture~\ref{conj_G Delta_UB}.
}\label{TableConstVal}
\end{table}

\section{Proof of the main result}
\label{sec:main}

In this section, we shall prove our main result, namely Theorem~\ref{thm:trianglefree}. Recall its statement.

\smallskip
\noindent \textbf{Theorem~\ref{thm:trianglefree}.}
Let $\Delta\ge 3$ be an integer, and let $G$ be a connected triangle-free graph of order $n \ge 3$. If $G \in \mathcal{F}_{\Delta}$, then
\[
\ID(G) = \left( \frac{\Delta - 1}{\Delta}\right) n+\frac{1}{\Delta}.
\]
On the other hand, if $G \notin \mathcal{F}_{\Delta}$ has maximum degree $\Delta$, then
\[
\ID(G) \le \left( \frac{\Delta - 1}{\Delta}\right) n.
\]

\noindent
\textbf{Proof of Theorem~\ref{thm:trianglefree}.} The first part of the statement holds by Proposition~\ref{prop-F_Delta}.

For the second part, let $G$ be a connected triangle-free graph of order $n$ and size~$m$ with maximum degree $\Delta \ge 3$ such that $G \notin \mathcal{F}_{\Delta}$. Thus, $n\geq5$. We proceed by induction on $n + m$ to show that $\ID(G) \le \left( \frac{\Delta - 1}{\Delta}\right) n$. Since $G$ is connected, we note that $m \ge n-1$, and so $n + m \ge 2n-1 \ge 9$. If $n + m = 9$, then $G$ is formed from a star by joining a pendant leaf to another leaf. By Theorem \ref{thm_trees}, we have $\ID(G)\le \frac{2}{3}n$, in this case. Furthermore, if $n =m=5$, then $G$ is formed from a $4$-cycle joined by a leaf. Observe that such a graph has $\ID(G)=3<\frac{2}{3}n$. This establishes the base cases. Let $n + m \ge 11$, where $n \ge 5$.

For the inductive hypothesis, assume that if  $G'$ is a connected triangle-free graph of order $n' \ge 3$ and size~$m'$ with $n' + m' < n + m$ and with   maximum degree $\Delta'=\Delta(G') \ge 3$  such that $G' \notin \mathcal{F}_{\Delta'}$, then $\ID(G') \le \left( \frac{\Delta' - 1}{\Delta'}\right) n'$.

If $m = n-1$, then $G$ is a tree. Since $G \notin \mathcal{F}_{\Delta}$, by Theorem~\ref{thm_trees}, we have  $\ID(G) \le \left( \frac{\Delta - 1}{\Delta}\right) n$. Hence, we may assume that $m \ge n$, for otherwise the desired upper bound follows. Thus, the graph $G$ contains a cycle edge, that is, an edge that belongs to a cycle in $G$. Moreover, since $\Delta \ge 3$, the graph $G$ is not a cycle.

Among all cycle edges in $G$, let $e = uv$ be chosen so that the sum of the degrees of its ends is as large as possible, that is, $\deg_G(u) + \deg_G(v)$ is maximal. Since $\Delta(G) \ge 3$, we have for the edge $e$ that $$\deg_G(u) + \deg_G(v) \ge 5.$$ Let $$G' = G - e.$$ Since $e$ is a cycle edge of $G$, the graph $G'$ is a connected triangle-free  graph of order~$n$. Let $\Delta(G') = \Delta'$, and so $\Delta' \ge \Delta - 1$. We note that $n(G') = n$ and $m(G') = m - 1$.

Suppose that $\Delta' = 2$. In this case, $\Delta = 3$ and $G'$ is either a path or a cycle. Suppose firstly that $G'$ is a cycle. By the triangle-free condition and our choice of the edge $e$, we infer that $G$ is obtained from a cycle $C_n$ where $n \ge 6$ by adding a chord between two non-consecutive vertices on the cycle in such a way as to create two cycles that both contain the edge~$e$ and both have length at least~$4$. Assume first that $n \neq 7$. By Corollary~\ref{cor_paths & cycles}(f), $\ID(G') \le \frac{2}{3}n$. By~Observation~\ref{ob:Ftree1},
$\ID(G)\le\ID(G')\le \frac{2}{3}n$ if $n$ is odd. Hence, we may consider the case with even $n \ge 6$. Let $G'$ be the cycle $v_1 v_2 \cdots v_n v_1$. Observe that both the set $V_{\even}$ of vertices with even subscript and the set $V_{\odd}$ of vertices with odd subscript are identifying codes of $G'$ of size $n/2$. If the chord in $G$ is between two vertices with even subscripts, then the set $V_{\odd}$ is an identifying code in $G$, and vice versa. Moreover, if the chord is between a vertex with even subscript and a vertex with odd subscript, then both sets $V_{\even}$ and $V_{\odd}$ are identifying codes of $G$. Hence, $\ID(G) \le \frac{1}{2}n$. If $n = 7$, then $\ID(G) = \ID(G') - 1 = 4 < \frac{2}{3}n$ by Observation~\ref{ob:Ftree1}. Therefore, if $G'$ is a cycle, then $\ID(G) \le \frac{2}{3}n$, as desired. Suppose secondly that $G'$ is a path. In this case, $n \ge 5$. For small values of $n$, namely $n \in \{5,6,7,8,10\}$, it can readily be checked that $\ID(G)\le \ID(G') = \lfloor \frac{n}{2} \rfloor + 1 \le \frac{2}{3}n$. Hence we may assume that $n \ge 9$ is odd or $n \ge 12$ is even. Any identifying code in the path $G'$ can be extended to an identifying code in $G = G' + e$ by adding at most one vertex, and so, by Corollary~\ref{cor_paths & cycles}(a) we have $\ID(G) \le \ID(G') + 1 = \lfloor \frac{n}{2} \rfloor + 2 \le \frac{2}{3}n$. Hence, we may assume that $\Delta' \ge 3$, for otherwise the desired bound holds.

Since $G$ is triangle-free, we note that $G'\ne K_{1,\Delta'}$ (since otherwise adding back the deleted edge $e$ would create a triangle in $G$). In particular if $\Delta' \ge 4$, then $G' \notin \mathcal{F}_{\Delta'}$. If $\Delta' = 3$ and $G' \in \mathcal{F}_{\Delta'}$, then by Observation~\ref{ob:Ftree1} and Lemma~\ref{ob:treeFam} we infer that $\ID(G) < \frac{2}{3}n \le \left( \frac{\Delta - 1}{\Delta}\right) n$. Hence we may assume that $G' \notin \mathcal{F}_{\Delta'}$, for otherwise the desired bound holds. Applying the inductive hypothesis to the graph $G'$, we have $\ID(G') \le \left( \frac{\Delta' - 1}{\Delta'}\right) n \le \left( \frac{\Delta - 1}{\Delta}\right) n$.

For notational convenience, let $N_u = N_G(u) \setminus \{v\}$ and let $N_v = N_G(v) \setminus \{u\}$. Since $G$ is triangle-free and $uv \in E(G)$, we note that $N_u \cap N_v = \emptyset$. Let $A$ be the boundary of the set $\{u,v\}$, that is, $A = N_u \cup N_v$ is the set of vertices different from $u$ and $v$ that are adjacent to $u$ or $v$. Further, let
\[
A_{uv} = A \cup \{u,v\} = N_G[u] \cup N_G[v] \hspace*{0.5cm} \mbox{and} \hspace*{0.5cm} \barA_{uv} = V(G) \setminus A_{uv}.
\]

\begin{figure}[!htb]
\centering
\centering
\begin{tikzpicture}[
blacknode/.style={circle, draw=black!, fill=black!, thick},
whitenode/.style={circle, draw=black!, fill=white!, thick},
scale=1]
\node[blacknode,label={90:$u$}] (u) at (0,0) {};
\node[blacknode,label={90:$v$}] (v) at (3,0) {};

\node[whitenode] (2) at (-2,-2) {};
\node (2a) at (-1,-2) {};
\node (2b) at (0,-2) {};
\node[whitenode] (3) at (1,-2) {};
\node at ($(2)+(-0.75,0)$) {$N_u$};
\node at ($(2)+(1.5,0)$) {$\cdots$};
\draw[rounded corners] ($(2)+(-1.25,0.5)$) rectangle ($(3)+(0.5,-0.5)$);

\node[whitenode] (4) at (2.5,-2) {};
\node (4a) at (3.5,-2) {};
\node[whitenode] (5) at (4.5,-2) {};
\node at ($(5)+(0.75,0)$) {$N_v$};
\node at ($(4)+(1,0)$) {$\cdots$};
\draw[rounded corners] ($(4)+(-0.5,0.5)$) rectangle ($(5)+(1.25,-0.5)$);

\draw[rounded corners] ($(2)+(-1.5,0.75)$) rectangle ($(5)+(2.25,-0.75)$);
\node at ($(5)+(1.75,0)$) {$A$};

\draw[rounded corners] ($(2)+(-1.75,2.75)$) rectangle ($(5)+(2.5,-1)$);
\node at ($(v)+(3.25,0.25)$) {$A_{uv}$};

\draw[rounded corners] ($(2)+(-1.75,-1.25)$) rectangle ($(5)+(2.5,-3)$);
\node at ($(5)+(1.75,-2.25)$) {$\barA_{uv}$};

\node[whitenode] (6a) at (-1,-4.25) {};
\node[whitenode] (6b) at ($(6a)+(4,0)$) {};
\node at ($(6a)+(2,0)$) {$\cdots$};

\draw[-, thick,dashed] (u) -- (v) node[midway,above] {$e$};
\draw[-, thick] (u) -- (2);
\draw[-, thick] (u) edge ($(u)!0.45!(2a)$) edge [dotted] ($(u)!0.55!(2a)$);
\draw[-, thick] (u) edge ($(u)!0.45!(2b)$) edge [dotted] ($(u)!0.55!(2b)$);
\draw[-, thick] (v) edge ($(v)!0.45!(4a)$) edge [dotted] ($(v)!0.55!(4a)$);
\draw[-, thick] (2) edge ($(2)!0.65!(6a)$) edge [dotted] ($(2)!0.75!(6a)$);
\draw[-, thick] (3) edge ($(3)!0.65!(6a)$) edge [dotted] ($(3)!0.75!(6a)$);
\draw[-, thick] (3) edge ($(3)!0.4!(4)$) edge [dotted] ($(3)!0.55!(4)$);
\draw[-, thick] (4) edge ($(4)!0.65!(6a)$) edge [dotted] ($(4)!0.75!(6a)$);
\draw[-, thick] (5) edge ($(5)!0.65!(6b)$) edge [dotted] ($(5)!0.75!(6b)$);
\draw[-, thick] (2) edge ($(2)!0.65!(6b)$) edge [dotted] ($(2)!0.75!(6b)$);
\draw[-, thick] (4) edge ($(4)!0.65!(6b)$) edge [dotted] ($(4)!0.75!(6b)$);
\draw[-, thick] (u) -- (3);
\draw[-, thick] (v) -- (4);
\draw[-, thick] (v) -- (5);
\end{tikzpicture}
\caption{The general setting of the proof of Theorem~\ref{thm:trianglefree}. Since $G$ is triangle-free, $N_u$ and $N_v$ are independent sets, but there can be edges across them.}
\label{fig:general_setting}
\end{figure}
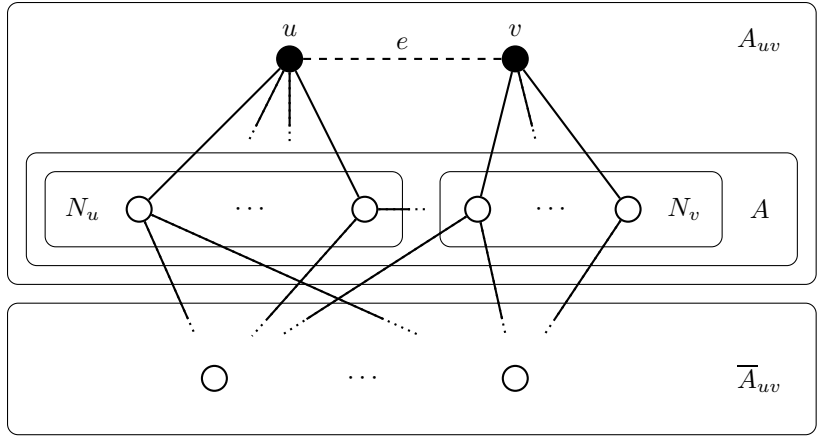

See Figure~\ref{fig:general_setting} for an illustration.

Let $C'$ be an optimal identifying code in $G'$, and so $C'$ is an identifying code in $G'$ and $|C'| = \ID(G')$. If $C'$ is an identifying code in $G$, then $\ID(G) \le |C'| = \ID(G') \le \left( \frac{\Delta - 1}{\Delta}\right) n$. Hence, we may assume that $C'$ is not an identifying code in $G$, for otherwise the desired bound holds.

We will next proceed with proving a series of claims.

\begin{claim}
\label{c:claimA}
If $V(G) = A_{uv}$, then
\[
\ID(G) \le \left( \frac{\Delta - 1}{\Delta}\right) n.
\]
\end{claim}
\proof Suppose that $V(G) = A_{uv}$. Recall that by our choice of the edge~$e$, we have $\deg_G(u) \ge 2$, $\deg_G(v) \ge 2$, and $5 \le \deg_G(u) + \deg_G(v) \le 2\Delta$. In this case since $V(G) = A_{uv}$, we have $n = \deg_G(u) + \deg_G(v) \le 2\Delta$. Let $u'$ be an arbitrary neighbor of $u$ in $G$ different from $v$, and let $v'$ be an arbitrary neighbor of $v$ in $G$ different from $u$. The code $C = V(G) \setminus \{u',v'\}$ is an identifying code in $G$, implying that
\[
\begin{array}{lcl}
\ID(G) \le |C| & = & \deg_G(u) + \deg_G(v) - 2 \2 \\
& = & \displaystyle{ \left( \frac{\deg_G(u) + \deg_G(v) - 2}{\deg_G(u) + \deg_G(v)} \right) n } \2 \\
& \le & \displaystyle{  \left( \frac{2\Delta - 2}{2\Delta} \right) n } \2 \\
& = & \displaystyle{ \left( \frac{\Delta - 1}{\Delta} \right) n, }
\end{array}
\]
yielding the desired upper bound.~\smallqed

\medskip

By Claim~\ref{c:claimA}, we may assume that $V(G) \ne A_{uv}$, for otherwise the desired result follows. Let
\[
G_{uv} = G - (N_G[u] \cup N_G[v]),
\]
and so $V(G_{uv}) = \barA_{uv} = V(G) \setminus (A \cup \{u,v\})$. Since $C'$ is an identifying code in $G'$ but not in $G$, a pair of vertices in $G$ is not identified by the code $C'$ when adding back the deleted edge $e$ to the graph $G'$ to reconstruct $G$. Since $G$ is triangle-free, the only possible pairs of vertices not identified by the code $C'$ in $G$ are the pairs $\{u,v\}$ or $\{u',v\}$ or $\{u,v'\}$ where $u'$ is some neighbor of $u$ different from $v$, and $v'$ is some neighbor of $v$ different from $u$.

\begin{claim}
\label{c:claimB}
If for every optimal identifying code $C'$ in $G'$ the pair $\{u,v\}$ is the only pair not identified by the code $C'$ in $G$, then
\[
\ID(G) \le \left( \frac{\Delta - 1}{\Delta}\right) n.
\]
\end{claim}
\proof Let $C'$ be an optimal identifying code in $G'$, and suppose the pair $\{u,v\}$ is the only pair not identified by the code $C'$ in $G$. Necessarily, $\{u,v\} \subseteq C'$, no neighbor of $u$ different from $v$ belongs to $C'$ and no neighbor of $v$ different from $u$ belongs to $C'$; that is, $C' \cap A = \emptyset$, and so no neighbor of $u$ in $G'$ and no neighbor of $v$ in $G'$ belongs to $C'$. Equivalently, $C' \setminus \{u,v\} = C'  \cap \barA_{uv}$. See Figure~\ref{fig:general_setting}, where the black vertices belong to $C'$, and the other ones do not. We proceed further with a series of subclaims that we will need when proving Claim~\ref{c:claimB}. Since every neighbor of $u$ (respectively, $v$) is identified by the code $C'$ in the graph $G'$, we infer the following claim.

\begin{subclaim}
\label{c:claimB.1}
Every neighbor of $u$ (respectively, $v$) in $G'$ has at least one neighbor that belongs to the set $C' \setminus \{u,v\}$.
\end{subclaim}

We shall frequently use the following claim when obtaining structural properties of the graph $G$.

\begin{subclaim}
\label{c:claimB.2}
If there exists a vertex $w \in C' \setminus \{u,v\}$ and a vertex $z \in A$ such that $(C' \setminus \{w\}) \cup \{z\}$ is an identifying code in the graph $G$, then $\ID(G) \le \left( \frac{\Delta - 1}{\Delta}\right) n$.
\end{subclaim}
\proof Consider the set $(C' \setminus \{w\}) \cup \{z\}$ be as defined in the statement of the claim. If this set is an identifying code in the graph $G$, then $\ID(G) \le |C'| \le \left( \frac{\Delta - 1}{\Delta}\right) n$.~\smallqed

\medskip

Thus, Claim~\ref{c:claimB.2} proves the result when the graph $G_{uv}$ contains a component of order~$1$. By Claim~\ref{c:claimB.2}, we may assume that there does not exist a vertex $w \in C' \setminus \{u,v\}$ and a vertex $z \in A$ such that $(C' \setminus \{w\}) \cup \{z\}$ is an identifying code in the graph $G$, for otherwise the desired bound holds. Recall that $C' \cap A = \emptyset$.  In the following claim, we show that there are no $P_2$-components in $G_{uv}$.

\begin{subclaim}
\label{c:claimB.3}
No component in $G_{uv}$ has order $2$.
\end{subclaim}
\proof Suppose that the graph $G_{uv}$ contains a component $F$ of order~$2$, and so $F$ is isomorphic to $P_2$. As observed earlier, in the graph $G'$ we have $C' \cap N(u) =  C' \cap N(v) = \emptyset$. In this case, the two vertices in the component $F$ are not separated by the code $C'$, a contradiction.~\smallqed

\medskip

In the following subclaims we consider the case with $\Delta=3$ separately. This is due to the more complex structure of set $\mathcal{F}_3$ compared to sets $\mathcal{F}_i$ for $i\geq4$.

\begin{subclaim}
\label{c:claimB.4}
Let $\Delta=3$. If $C_7$ is a component in $G_{uv}$, then $\ID(G) \le \left( \frac{\Delta - 1}{\Delta}\right) n$.
\end{subclaim}
\proof Let $F=C_7$ be a cycle component in $G_{uv}$.  Since $C'$ is an optimal identifying code in $G'$ and since $C'$ contains no vertices in the boundary $A$, the set $V(F) \cap C'$ is an  identifying code in $F$, that is, $\ID(F) \le |V(F) \cap C'|$. By Theorem~\ref{thm_GMS2006}, we have $\ID(C_7)=5$. Let $V(F)=\{w_1,w_2,\dots,w_7\}$ and $E(F)=\{w_7w_1\}\cup\{w_iw_{i+1}\mid 1\le i\le 6\}$. Observe that each vertex in $F$ can be adjacent to at most one vertex in $A$ since $\Delta=3$. Note that, since $\Delta=3$, we have $1\leq|N_u|\le 2$ and $1\leq|N_v|\le 2$, and so $|A|\le 4$. We denote $N_u=\{u_1,u_2\}$ and $N_v=\{v_1,v_2\}$ (if these vertices exist). We further assume that $|N(F)\cap N_u|\ge |N(F)\cap N_v|$. Moreover, let $w_2$ be adjacent to $u_1\in N_u$.   Let us consider vertex sets $$C_1=(C'\setminus V(F))\cup \{u_1,w_1,w_4,w_5,w_6\}$$ and $$C_2=(C'\setminus V(F))\cup \{u_1,w_3,w_5,w_6,w_7\}.$$ In the following, we show that at least one of these two sets is an identifying code in $G$. Notice that $|C_1|=|C_2|\le |C'|$.

If $u_1$ is the only vertex of $A$ adjacent to a vertex in $F$, then $C_1$ is an identifying code in $G$. In particular, we may use Lemma~\ref{LemP3} to see that $u,v$ and $u_1$ have unique neighborhoods in $C_1$, while $C_1\cap V(F)$ forms an identifying code for $F\setminus \{w_2\}$, and $w_2$ is identified by vertices $u_1$ and $w_1$. The remaining vertices are identified by the vertices in $C'$.

Assume then that there are two vertices of $A$ ($u_1$ and, say, $x$) adjacent to vertices in $F$. If $x$ is in $N_u$ ($x=u_2$), then $C_1$ is an identifying code in $G$, even if $x$ is not dominated by a vertex in $C_1\cap V(F)$. Indeed, we have $u\in N[u_2]\cap C_1$ but $v,u_1\not\in N[u_2]\cap C_1$. Thus, $u_2$ is separated from all other vertices. If $x\in N_v$, then $C_1$ or $C_2$ is an identifying code in $G$. In this case, we choose such a set $C_i$ ($i\in\{1,2\}$) so that $N(x)\cap V(F)\cap C_i\neq\emptyset$.

Assume next that there are three vertices of $A$ adjacent to vertices in $F$. In this case, due to our assumption that $|N_u\cap N(F)|\ge |N_v\cap N(F)|$, 
we have $|N_u\cap N(F)|=2$ 
and hence, we may assume that $A\cap N(F)=\{u_1,u_2,v_1\}$. Again we choose a set $C_i$ ($i\in\{1,2\}$), such that $N(v_1)\cap V(F)\cap C_i\neq\emptyset$, as our identifying code. Note that as in the previous case, we do not need to dominate $u_2$ from $F$.

Finally, when (all) four vertices of $A$ are adjacent to vertices in $F$, we again choose a set $C_i$ ($i\in\{1,2\}$) such that $N(v_1)\cap V(F)\cap C_i\neq\emptyset$ as our identifying code. As in the previous cases, we do not need to dominate vertices $u_2$ and $v_2$ from $F$. The argument for $u_2$ is similar as in the previous cases. Moreover, by Lemma~\ref{LemP3}, vertices $v$ and $u$ have unique neighborhoods in $C_i$, vertex $v$ separates $v_2$ from vertices other than $u,v$ and $v_1$, while $v_1$ is separated from $v_2$ by a vertex in $F$.

The claim follows in these cases since $|C_1|=|C_2|\le |C'|$ and since $C_1$ or $C_2$ is an identifying code in $G$.~\smallqed

\medskip

\begin{subclaim}
\label{c:claimB.5}
Let $\Delta=3$. If $C_4$ is a component in $G_{uv}$, then $\ID(G) \le \left( \frac{\Delta - 1}{\Delta}\right) n$.
\end{subclaim}
\proof Let $F=C_4$ be a cycle component in $G_{uv}$.  Since $C'$ is an optimal identifying code in $G'$ and since $C'$ contains no vertices in the boundary $A$, the set $V(F) \cap C'$ is an identifying code in $F$, that is (by Theorem~\ref{thm_GMS2006}), $\ID(F)=3 \le |V(F) \cap C'|$. Let $V(F)=\{w_1,w_2,w_3,w_4\}$ and $E(F)=\{w_1w_2,w_2w_3,w_3w_4,w_4w_1\}$. Observe that each vertex in $F$ can be adjacent to at most one vertex in $A$ since $\Delta=3$. Recall that by Claim~\ref{c:claimB.1} every vertex in $A$ is dominated by two vertices in $C'$.

Let us first assume that $V(F)\subseteq C'$. Then, we can apply Claim~\ref{c:claimB.2} to any vertex in $F$ and its neighbor in $A$ and the claim follows. Thus, we can now assume that $|V(F)\cap C'|=3$.

Let us first assume that there are one or three vertices in $A$ adjacent to vertices in $F$. Note that each vertex in $F$ can be adjacent to at most one vertex in $A$. If there is one vertex in $A$ adjacent to a vertex in $F$, then we assume that the adjacent vertex in $F$ is $w_2\in F$. If there are three vertices in $A$ adjacent to vertices in $F$, then we assume that they are adjacent to vertices $w_1,w_2$ and $w_3$ in $F$. Furthermore, we assume that $z\in A$ is adjacent to $w_2$. We will consider vertex set $C=(C'\setminus V(F))\cup\{w_1,w_3,z\}$. This is an identifying code in $G$ since all vertices in $A_{uv}$ are separated from other vertices by $u$ and $v$ and vertices $u,v$ and $z$ have unique code neighborhoods in $C$ by Lemma~\ref{LemP3}. Moreover, vertices in $F$ are separated from each other by $z,w_1$ and $w_3$ and finally vertices in $A\setminus\{z\}$ are separated from each other by the same vertices as in $C'$.

Assume then that there are four vertices in $A$ adjacent to vertices in $F$. Let us assume without loss of generality that $\{w_1,w_2,w_3\}\subset C'$ and that $z\in A$ is adjacent to $w_2$. As in the previous case, now $C=(C'\setminus\{w_2\})\cup\{z\}$ is an identifying code in $G$.

Finally, we have the case where we have exactly two vertices in $A$ adjacent to vertices in $F$. Assume first that these vertices in $A$ are adjacent to two adjacent vertices of $F$, say to $w_1$ and $w_2$. Let $z\in A$ be a vertex adjacent to $w_2$. As in the previous cases, we let $C=(C'\setminus V(F))\cup \{w_1,z,w_3\}$. With the same arguments as in the previous cases, $C$ is an identifying code in $G$.

Let us next consider the case where there are two vertices in $A$ adjacent to only non-adjacent vertices in $F$, say vertices $w_1$ and $w_3$. Observe that now $w_2$ and $w_4$ are open twins in $G'$ and we have exactly two edges between $A$ and $F$. Let $z\in A$ be a neighbor of $w_1$. Assume that $|N_u|\geq|N_v|$. Let us denote the neighbors of $u$ in $A$ by $u_1$ and $u_2$ and neighbors of $v$ in $A$ by $v_1$ and $v_2$ (if $v_2$ exists). Notice that if both edges from $F$ to $A$ are between $F$ and $N_u$ or $F$ and $N_v$, then $C=\{z,w_2,w_4\}\cup (C'\setminus V(F))$ is an identifying code of $G$. Thus, we may assume without loss of generality, that there is an edge from $w_1$ to $u_1$ and from $w_3$ to $v_1$. Assume next that for   $u_2$ or $v_2$ there exists vertex $c\in C'$ such that $c\in N(v_2)\setminus N(v_1)$ or $c\in N(u_2)\setminus N(u_1)$. If the former holds, then $C_u=\{u_1,w_2,w_4\}\cup (C'\setminus V(F))$ is an identifying code of $G$ and if the latter holds, then $C_v=\{v_1,w_2,w_4\}\cup (C'\setminus V(F))$ is an identifying code of $G$. Hence, we next assume that $c$ does not exist. Moreover, if $v_2$ does not exist, then again $C_u$ is an identifying code of $G$ and we are done. Hence, we assume that $v_2$ exists. Recall that each vertex in $A$ is dominated by at least two vertices of $C'$ (since otherwise, a vertex of $A$ dominated only by one vertex of $C'$ would not be separated from either $u$ or $v$ in $G'$). Hence, $u_2$ has a neighbor in $C'\cap \overline{A}_{uv}$, and since every such vertex is also a neighbor of $u_1$ (by the above argument that $c$ does not exist), there exists a vertex $c_u\in C'\cap (N(u_1)\cap N(u_2))\setminus \{u\}$. By the same argument applied to $v_2$, there is also a vertex $c_v\in C'\cap (N(v_1)\cap N(v_2))\setminus \{v\}$. Observe that if there exists a vertex $w_u$ such that $N[w_u]\cap C'=\{c_u\}$ or a vertex $w_v$ such that $N[w_v]\cap C'=\{c_v\}$, then $w_u=c_u$ and $w_v=c_v$. Indeed, otherwise we would have  $N[w_u]\cap C'=N[c_u]\cap C'$ or $N[w_v]\cap C'=N[c_v]\cap C'$. We may also observe that we have either $\{w_1,w_2,w_3\}\subseteq C'$ or $\{w_1,w_4,w_3\}\subseteq C'$. Let us assume, without loss of generality, the former. However, now $C=\{u_1\}\cup (C'\setminus \{u\})$ is an identifying code in $G$. Indeed, only $u$, $u_1$ and $u_2$ have lost a codevertex from their neighborhoods compared to $C'$. Moreover, $u$ is the only vertex adjacent to both $v$ and $u_1$, $u_1$ is the only vertex adjacent to $w_1$ and $c_u$ and finally, $C\cap N[u_2]=\{c_u\}$ but now $C\cap N[c_u]\supseteq\{c_u,u_1\}$. Thus, also $u_2$ has a unique code neighborhood. Now, the claim follows.~\smallqed

\medskip

\begin{subclaim}
\label{c:claimB.6}
Let $\Delta=3$. If $P_4$ is a component in $G_{uv}$, then $\ID(G) \le \left( \frac{\Delta - 1}{\Delta}\right) n$.
\end{subclaim}
\proof Let $F=P_4$ be a path component in $G_{uv}$.  Since $C'$ is an optimal identifying code in $G'$ and since $C'$ contains no vertices in the boundary $A$, the set $V(F) \cap C'$ is an identifying code in $F$, that is, by Theorem~\ref{thm_BCHL2004}, $\ID(F)=3 \le |V(F) \cap C'|$. Let $V(F)=\{w_1,w_2,w_3,w_4\}$ and $E(F)=\{w_1w_2,w_2w_3,w_3w_4\}$.

\smallskip

\noindent\emph{Case 1: there is an edge between $w_1$ or $w_4$ and $A$.} Without loss of generality, let it be $w_1$. Consider now the graph $G_1=G-\{w_2,w_3,w_4\}$. First observe that $G_1$ is connected. Second, if $G_1\not\in \mathcal{F}_3$, then by induction we have an identifying code $C_1\subseteq V(G_1)$ of cardinality at most $\frac{2}{3} n(G_1)$. Moreover, if no neighbor (in $G$) of $w_2,w_3$ or $w_4$ is in $C_1$, then $C_1\cup\{w_2,w_4\}$ is an identifying code of claimed cardinality in $G$. If on the other hand $C_1\cap(N_G(w_2)\cup N_G(w_3))\neq\emptyset$, then at least one of the sets $C_1\cup\{w_2,w_4\}$ and $C_1\cup\{w_3,w_4\}$ is an identifying code in $G$. The case with $C_1\cap(N_G(w_4)\cup N(w_3))\neq\emptyset$ is analogous. Moreover, these identifying codes have the claimed cardinality.

Let us next consider the case with $G_1\in\mathcal{F}_3$. Furthermore, there is a vertex with degree~$3$ and there are at least six vertices in $G_1$ since $\deg(u)+\deg(v)\geq5$ by our assumption. Thus, $G_1\in \mathcal{F}_3 \setminus \{P_4,C_4,C_7,K_{1,3}\}=\mathcal{T}_3 \setminus \{K_{1,3}\}$. Let us assume first that  $G_1=T_2$. In this case, we have $n=10$.
Notice that $\deg(u)+\deg(v)=5$ and either $u$ or $v$ is the support vertex of degree $3$ in $T_2$ and the other one is the adjacent vertex of degree $2$. However, since $w_1$ is a leaf in $T_2$, this means that it must have distance of either~1 or~3 to set $\{u,v\}$, a contradiction since that distance is actually~2.

Let us next assume that $G_1=T_3$. Observe that in this case $\deg(u)+\deg(v)=6$. Since $w_1$ has distance~$2$ to set $\{u,v\}$, vertex $w_1$ has to be one of the leaves adjacent to one of the support vertices of degree~3 of $T_3$. But then, the leaf of $T_3$ adjacent to the unique support vertex of degree~2 in $T_3$ as well as that support vertex form a component of order~2 in $G_{uv}$, contradicting Claim~\ref{c:claimB.3}.

Assume next that $G_1\in\mathcal{T}_3 \setminus\{K_{1,3},T_2,T_3\}$. Let $C_1$ be an optimal identifying code of $G_1$ such that every vertex with degree at most~2 in $G_1$ is in $C_1$ (such a code exists by Proposition~\ref{propTreeCodes}, see Figure~\ref{fig:trees}).
 Since $G$ is not a tree, there are at least two edges between  $A\cup\{w_1\}$ and $\{w_2,w_3,w_4\}$.
Hence,  there is an edge between $w_1$ and $w_2$ and between $G_1'=G_1-\{w_1\}$ and $\{w_2,w_3,w_4\}$. 

Let us consider the case where there are no edges between $A$ and $w_2$ or $w_3$ but there is an edge to $w_4$. In this case, we may consider the identifying code $C=C_1\cup\{w_3\}$ which has the claimed cardinality. Notice that this is indeed an identifying code since $w_1\in C_1$ and $w_4$ is dominated by some vertex in $C_1$.	

Therefore, we may assume from now on that there is an edge between $A$ and $w_2$ or $w_3$. Let us assume that the edge is from $A$ to $w_2$ (the case with an edge to $w_3$ is similar). Consider now graph $G_1''=G-\{w_1,w_3,w_4\}$ together with an optimal identifying code $C_1''$. Observe that $G_1''$ is a tree since $G_1$ is a tree and $w_2$ has only one edge to $A$. Moreover, we also notice that there are no edges between vertices in $A$ since $G_1$ is a tree. Thus, if a 4-cycle contains exactly two vertices in $F$, then those vertices are $w_1$ and $w_4$.

If $G_1''\not\in \mathcal{F}_3$, then we have two cases based on whether $w_2\in C_1''$.

If $w_2\in C_1''$, then $C_1''\cup\{w_3,w_4\}$ is an identifying code in $G$. Indeed, by Lemma~\ref{LemP3} vertices $w_2,w_3$ and $w_4$ have unique code neighborhoods. Moreover, also $w_1$ has a unique neighborhood in $C_1''\cup\{w_2\}$ since any vertex of $A$ adjacent to $w_2$ is also either in $C_1''$ or has another neighbor outside of $A$ in $C_1''$.

If $w_2\not\in C_1''$, then we use set $C_{23}=C_1''\cup\{w_2,w_3\}$. Since $C_1''$ is an identifying code in $G_1''$, the vertex adjacent to $w_2$ (call it $z$) in $A$ is in $C_1''$ and $z$ has another code neighbor. Hence, by Lemma~\ref{LemP3}, vertices $z,w_2$ and $w_3$ have unique code neighborhoods. Moreover, since $w_1$ is adjacent to $w_2$ and all other vertices in $N[w_2]$ have unique code neighborhoods, also $w_1$ has a unique neighborhood in $C_{23}$. Furthermore, also $w_4$ has a unique neighborhood in $C_{23}$ since the only other vertex $x$ adjacent to $w_3$ which might not be separated from $w_4$ is in $A$ and is in $C_1''$ or has another vertex in $C_1''$ adjacent to it. Since $w_3$ and $w_4$ cannot belong to the same 4-cycle in $G$ and there are no triangles, vertices $w_4$ and $x$ are separated. All the other vertices are pairwise separated by the set $C_1''$. Hence, $C_{23}$ is an identifying code in $G$.

Thus, we may assume that $G_1''\in \mathcal{F}_3$. Observe that we may construct $G_1''$ from $G_1$ by removing a leaf and then adding a new leaf and vice versa. Hence, $G_1,G_1''$ are two trees of $\mathcal{F}_3$ (and thus, $\mathcal{T}_3$) with the same order of at least~6. Observe that for every $i$ such that $i\in\{0,1\}$ or $i\ge 4$, there are no support vertices of degree~2 in $T_i$ ($T_i\in\mathcal{T}_3$). However, when we remove a leaf from any such tree of $\mathcal{T}_3$ and add a different leaf to a vertex of degree at most~$2$, there necessarily exists a support vertex of degree~$2$ in the resulting tree. Thus, at least one of $G_1$ and $G_1''$ must be in $\{T_2,T_3\}$. Note that we cannot obtain an isomorphic copy of $T_3$ by this operation when starting from $T_3$, so at most one of $G_1$ and $G_1''$ is $T_3$. Thus, $\{G_1,G_1''\}=\{T_1,T_2\}$, $\{G_1,G_1''\}=\{T_2,T_2\}$ or $\{G_1,G_1''\}=\{T_3,T_4\}$. 
Let us first assume that $G_1,G_1''=\{T_3,T_4\}$. 
Consider the leaf adjacent to the support vertex of degree $2$ in $T_3$. Notice that the leaf has to belong to $F$. Moreover, its distance from set $\{u,v\}$ is exactly~2. However, now when we form $T_4$ by removing this leaf and attaching the new leaf, the new leaf is attached adjacent to $u$ or $v$ and hence, $A\cap V(F)\neq \emptyset$, a contradiction. Thus, $G_1,G_1''\in\{T_1,T_2\}$.

 We thus have $\{G_1,G_1''\}=\{T_1,T_2\}$ or $\{G_1,G_1''\}=\{T_2,T_2\}$. In both cases, we have $\deg(u)+\deg(v)=5$. Hence, one of $u$ or $v$ is the degree~$3$ support vertex in $T_2$ and the other one is the adjacent degree~$2$ vertex. Moreover, the leaf adjacent to the support vertex of degree~$2$ in $T_2$ is in $F$. But this vertex should have distance~2 to set $\{u,v\}$, a contradiction. This finishes the proof of Case~1.

\smallskip

\noindent\emph{Case 2: there are no edges from $w_1$ or $w_4$ to $A$.} Then, there is an edge from $A$ to $w_2$ or $w_3$; without loss of generality, assume there is an edge from $A$ to $w_2$ and denote $G_2=G-\{w_1,w_3,w_4\}$. Observe that if $G_2\not\in \mathcal{F}_3$, then by induction, there exists an identifying code $C_2$ in $G_2$ with $|C_2|\le \frac{2}{3} n(G_2)$ and $|C_2|+2\le \frac{2}{3} n$. Moreover, either $C_2\cup\{w_1,w_3\}$ or $C_2\cup\{w_2,w_3\}$ is an identifying code in $G$ depending on whether $w_2\in C_2$ or not. Hence, we may assume that $G_2\in \mathcal{F}_3$ and in fact $G_2\in \mathcal{T}_3$ since $G_2$ has a vertex of degree~3. If we do not have an edge from $w_3$ to $A$, then $G$ would be a tree, a contradiction. Hence, we may assume that there exists an edge from $w_3$ to $A$. Let us assume that $G_2\neq T_2$ and $C_2$ contains all vertices of degree at most~2 in $G_2$. This is possible by Proposition~\ref{propTreeCodes}. Observe that in this case $C_2\cup\{w_3\}$ is an identifying code of claimed size in $G$. Hence, we may assume that $G_2=T_2$. Since $\deg(u)+\deg(v)\geq5$, the only leaf of $T_2$ not in $A$ is the leaf adjacent to the degree~2 support vertex. However, this vertex has distance~3 to set $\{u,v\}$. Thus $G_2\neq T_2$, and the claim follows.~\smallqed

\medskip

\begin{subclaim}
\label{c:claimB.7}
Let $\Delta=3$. If $T\in \mathcal{T}_3$ is a component in $G_{uv}$, then $\ID(G) \le \left( \frac{\Delta - 1}{\Delta}\right) n$.
\end{subclaim}
\proof Let $F=T\in\mathcal{T}_3$ be a tree component in $G_{uv}$.  Since $C'$ is an optimal identifying code in $G'$ and since $C'$ contains no vertices in the boundary $A$, the set $V(F) \cap C'$ is an identifying code in $F$, that is, $\ID(F) \le |V(F) \cap C'|$. Let us call by $T^2$ the set of vertices of degree at most~2 in $G[F]$. Note that only vertices of $T^2$ can have a neighbor in $A$.

Let us first assume that $F\in \mathcal{T}_3\setminus\{T_2,T_3\}$.
By Proposition~\ref{propTreeCodes} we can choose an optimal identifying code $C_T$ for $T$ such that $T^2\subseteq C_T$ (see Figure~\ref{fig:trees}). Moreover, by Proposition~\ref{propTreeCodes}(ii) $C_T\setminus\{t\}$ is an optimal identifying code in $T-\{t\}$ for any $t\in T^2$. Notice that $C_T'=(C'\setminus V(F))\cup C_T$ is an identifying code in $G'$ and $|C'|=|C_T'|$. Notice that Claim~\ref{c:claimB.2} holds also for $C_T'$.
Assume first that $z\in A$ is adjacent to $t\in T^2$. If $C=\{z\}\cup C_T'\setminus\{t\}$ is an identifying code in $G$, then the claim follows from Claim~\ref{c:claimB.2}. If $C$ is not an identifying code, then we will show next that vertex $t$ must be adjacent to a vertex $u_1\in N_u$ and $v_1\in N_v$ and $t$ is the only vertex in $C_T'$ separating $u_1$ from another vertex $u_2\in N_u$, and similarly, vertex $v_2\in N_v$ from $v_1\in N_v$.
Indeed, since $u,v\in C$, all vertices of $A$ are separated from vertices in $F$. Moreover, since $C_T\setminus\{t\}$ is an identifying code in $F-\{t\}$, set $C\cap V(F)$ is an identifying code in $F-\{t\}$. Furthermore,  vertex $t$ is the only vertex which has exactly $z$ in its code neighborhood. Consequently, $z$, $u$ and $v$ are identified by Lemma~\ref{LemP3}. Therefore, the only vertices which might not be separated belong to $A$. Since $C_T'$ separated all vertices in $A$, we require $t$ to separate some vertices which are not separated by $z$. Moreover, $t$ can have at most two neighbors in $A$. If $t$ is adjacent to only one vertex in $A$, then that vertex is $z$ and $z\in C$ separates itself from other vertices in $A$. Thus, $t$ is adjacent to two vertices in $A$. If both of these vertices are in $N_u$ (or $N_v$), then $u$ and $v$ separate them from other vertices in $A$. Moreover, $z\in C$ separates it itself from other neighbors of $t$. Hence, $t$ is adjacent to a vertex $u_1\in N_u$ and $v_1\in N_v$. Since $u,v\in C$, the only vertex $u_1$ that can have the same code neighborhood with respect to $C$ is vertex  $u_2\in N_u$. A similar statement holds for $v_1$ and $v_2\in N_v$. 
Assume first that $C\setminus \{z\}$ separates $v_1$ and $v_2$ in $G$. In this case $(C\cup\{u_1\})\setminus \{z\}$ is an identifying code in $G$. A similar argument holds for $C\setminus \{z\}$ separating $u_1$ and $u_2$ in $G$. Hence, we may assume that $C\setminus \{z\}$ does not separate pairs $u_1,u_2$ and $v_1,v_2$. Thus, we require $t$ to separate these two pairs in $C_T'$, as claimed.

Since $v_2$ and $u_2$ are dominated by a vertex in $C_T'\setminus \{u,v\}$ (Claim~\ref{c:claimB.1}), there is a code vertex $w_u\in (C_T'\setminus\{u,v\})\cap N(u_1)\cap N(u_2)$ and a code vertex $w_v\in (C_T'\setminus\{u,v\})\cap N(v_1)\cap N(v_2)$. In this case, we may consider graph $G'$ and modify the identifying code $C_T'$ into $C''=\{u_1,v_1\}\cup C_T'\setminus\{u,v\}$ which is an identifying code in $G'$. This is a contradiction because $C''$ is an optimal identifying code of $G'$ such that $u,v$ are identified by $C''$ in $G$, contradicting the hypothesis of Claim~\ref{c:claimB}.

Assume next that $F=T_2$. Denote its unique support vertex of degree~3 by $s$. Let the leaves adjacent to $s$ be denoted by $l_1$ and $l_2$ and the third one be $l_3$. Finally, let $f_1,f_2,f_3$ be the three vertices on the path from $s$ to $l_3$ (in that order). Assume first that there are no edges between $A$ and $f_2,f_3$ or $l_3$. Now consider the graph $G_f=G-\{f_2,f_3,l_3\}$. Assume first that $G_f\not\in \mathcal{F}_3$ and let $C_f$ be an optimal identifying code in $G_f$. In this case, if $f_1\in C_f$, then $C_f\cup\{f_2,f_3\}$ is an identifying code of claimed cardinality in $G$. If $f_1\not\in C_f$, then $C_f\cup\{f_2,l_3\}$ is an identifying code in $G$ of claimed cardinality. Assume then that $G_f\in \mathcal{F}_3$. Since there is a vertex of degree~3 in $G_f$, we have $G_f\in \mathcal{T}_3$. However, since there are no edges from $A$ to $\{f_2,f_3,l_3\}$, graph $G$ is a tree, a contradiction. Hence, we may assume that there is an edge from $A$ to $\{f_2,f_3,l_3\}$.

Let us next consider graph $G_s=G-\{l_1,l_2,s\}$ and its optimal identifying code by $C_s$. Notice that $G_s$ is connected. Assume first that $G_s\not\in \mathcal{F}_3$. If one of the three vertices in $\{l_1,l_2,s\}$ is dominated by a vertex in $C_s$ in $G$, then the set $C_s$ together with two adjacent vertices from $\{l_1,l_2,s\}$ such that at least one of them is dominated by a vertex from $C_s$ is an identifying code of claimed cardinality. If no vertex in $\{l_1,l_2,s\}$ is dominated by a vertex from $C_s$, then we consider $C_s\cup\{l_1,l_2\}$. This is an identifying code, since $s$ is separated from all other vertices. Indeed, if some vertex in $V(G)\setminus V(F)$ is adjacent to $l_1$ and $l_2$, then it is adjacent also to some vertex in $C_s\cap (V(G)\setminus V(F))$ while vertex $s$ cannot be (since it is not dominated by $C_s$). Hence, we may assume that $G_s\in \mathcal{F}_3$ and more specifically $G_s\in \mathcal{T}_3$ since there is a vertex of degree~3. Hence, there is a single edge between $\{f_1,f_2,f_3,l_3\}$ and $A$. Furthermore, if there is an edge from $A$ to $f_2$ or $f_3$, then $G_s\not\in \mathcal{T}_3$, and if the edge is from $A$ to $f_1$ or $l_3$, then $G_s$ has to be $T_2$. However, there are at least five vertices in $A_{uv}$ and hence $G_s\neq T_2$. Thus $G_s\not\in \mathcal{F}_3$.

Let us finally consider the case where $F=T_3$. Denote by $s_1$ and $s_2$ its two support vertices of degree~3, and by $f_2$ the support vertex of degree~$2$. Let leaves $l_1$ and $l_2$ be adjacent to $s_1$, and leaves $l_3$ and $l_4$ be adjacent to $s_2$. Furthermore, let the leaf adjacent to $f_2$ be $l_5$ and the other vertex adjacent to $f_2$ be $f_1$. Further denote the vertex of degree~3 adjacent to $f_1$ by $f_s$. Assume first that there are no edges from $A$ to $\{f_1,f_2,l_5\}$. Thus, $G_f=G-\{f_1,f_2,l_5\}$ is a connected graph and let $C_f$ be an optimal identifying code in $G_f$. Notice that if $G_f\in \mathcal{F}_3$, then $G_f\in \mathcal{T}_3$, and then $G$ is a tree, a contradiction. Thus, $G_f\not\in \mathcal{F}_3$ and $C_f$ contains at most two-thirds of the vertices of $G_f$. In this case, if $f_s\in C_f$, then the set $C=C_f\cup\{f_1,f_2\}$ is an identifying code in $G$ and if $f_s\not\in C_f$, then the set $C=C_f\cup\{f_1,l_5\}$ is an identifying code in $G$. Moreover, both of these sets contain at most two-thirds of the vertices in $G$, and we are done. Hence, we may assume from now on that there is an edge from $A$ to $\{f_1,f_2,l_5\}$. Assume then that there are no edges from $A$ to $\{s_1,l_1,l_2\}$. Now, graph $G_s=G-\{s_1,l_1,l_2\}$ is connected, has a cycle and hence, by induction, also an optimal identifying code $C_s$ satisfying the two-thirds upper bound. Moreover, set $C_s\cup\{l_1,l_2\}$ is an identifying code in $G$.
Hence, there is an edge from $A$ to $\{s_1,l_1,l_2\}$. By symmetry, a similar argument holds for $\{s_2,l_3,l_4\}$, so there is an edge from $A$ to $\{s_2,l_3,l_4\}$.
Hence we can assume next that that there is an edge from $A$ to sets  $\{s_1,l_1,l_2\}$, $\{s_2,l_3,l_4\}$ and $\{l_5,f_1,f_2\}$. Let us consider graph $G_f'=G-\{s_2,l_3,l_4\}$ together with an optimal identifying code $C_f'$ of $G_f'$. Notice that $G_f'$ is connected, has a cycle, and maximum degree~3. Hence, $G_f'\not\in \mathcal{F}_3$ and by induction, $|C_f'|\le \frac{2}{3}n(G_f')$. Consider set $C_f'\cup\{l_3,l_4\}$ if no vertex in $C_f'$ dominates a vertex in $\{l_3,l_4\}$ and otherwise, set $C_f'$ together with $s_2$ and a vertex in $\{l_3,l_4\}$ that is dominated by $C_f'$. Notice that in each case, the corresponding set is an identifying code of claimed cardinality for $G$.
~\smallqed

\medskip

By the above claims, we may assume that if $\Delta=3$, then for any component $F$ in $G_{uv}$ we have $\ID(F)\le  \frac{2}{3} n(F)$.

\begin{subclaim}
\label{c:claimB.8}
Let $\Delta\geq4$. If $T$ is a $\Delta$-star component in $G_{uv}$, then $\ID(G) \le \left( \frac{\Delta - 1}{\Delta}\right) n$.
\end{subclaim}
\proof Let $T$ be a $\Delta$-star component in $G_{uv}$. Let $V(T)=\{w,w_1,\dots,w_\Delta\}$ and $w$ be the center vertex of $T$. Observe that $w$ is adjacent only to vertices in $T$. Let $w_1$ be adjacent to a vertex $z\in A$. Let $T'=T-\{w_1\}$ and let $G_T=G-T'$. Observe that $G_T$ is a connected graph and at least one of $u$ or $v$ has degree at least~$3$ in $G_T$. Hence, there are at least six vertices in $G_T$. Let us assume first that $G_T\in \mathcal{F}_3$; then, $G_T\in \mathcal{T}_3$. In this case we have $\ID(G_T)\le \frac{3}{4}n(G_T)$ and $\ID(T')\leq\left(\frac{\Delta-1}{\Delta}\right)n(T')$. Moreover, there are at least three leaves in $T'$. Let $C_T$ be an optimal identifying code in $G_T$. Observe that $C=C_T\cup\{w,w_2,\dots,w_{\Delta-1}\}$ is an identifying code in $G$ of cardinality at most $|C|\le \frac{3}{4}n(G_T)+ \left(\frac{\Delta-1}{\Delta}\right)n(T')\le \left(\frac{\Delta-1}{\Delta}\right)n$. Furthermore, if $G_T\not\in \mathcal{F}_3$, since also
$G_T \notin \mathcal{F}_{\Delta'}$, where $\Delta' = \Delta(G_T)$, there exists an identifying code $C_T$ of $G_T$ such that $|C_T|=\ID(G_T)\leq \left(\frac{\Delta'-1}{\Delta'}\right)n(G_T) \leq\left(\frac{\Delta-1}{\Delta}\right)n(G_T)$ and hence, we may again consider $C=C_T\cup\{w,w_2,\dots,w_{\Delta-1}\}$ as our identifying code for $G$ and $|C|\le \left(\frac{\Delta-1}{\Delta}\right)n$ which completes the proof of the claim.~\smallqed  \medskip

Let us denote the set of isolated vertices in $G_{uv}$ by $\mathcal{I}$ and let us further denote $G_{AI}=G[A_{uv}\cup\mathcal{I}]$.

\begin{subclaim}
\label{c:claimB.9}
Set $A_{uv}$ is an identifying code in $G_{AI}$.
\end{subclaim}
\proof Observe that $\mathcal{I}\subseteq C'$ (otherwise some vertex of $\mathcal I$ would not be dominated by $C'$ in $G'$). If set $A_{uv}$ is not an identifying code in $G_{AI}$, then there are open twins $w_1,w_2$ in $\mathcal{I}$. 
Let $z\in N(w_1)$. Notice that $z\in A$ and hence, $z\not\in C'$. Since vertices in $\{w_1,w_2\}$ are open twins, we may consider set $C=\{z\}\cup C'\setminus \{w_1\}$. Observe that $C$ is an identifying code of $G$ since $w_2$ separates all the same vertices as vertex $w_1$.
Hence, we obtain a contradiction from Claim~\ref{c:claimB.2}. Thus, there are no twins in $\mathcal{I}$ and set $A_{uv}$ is an identifying code in $G_{AI}$.~\smallqed

\medskip

 Let us next consider the minimum cardinality of an identifying code in $G_{AI}$.
 Assume that $|N_u|\ge |N_v|$. Since $|N_v|\leq|N_u|\leq\Delta-1$, we have $|A_{uv}|\le 2\Delta$. Thus, if $|A_{uv}|>\left( \frac{\Delta - 1}{\Delta}\right)n(G_{AI})=\left( \frac{\Delta - 1}{\Delta}\right)(|A_{uv}|+|\mathcal{I}|)$, then $|A_{uv}|>(\Delta-1)|\mathcal{I}|$. Hence,  this implies that $|\mathcal{I}|\leq2$. By Claim~\ref{c:claimB.9}, the set $\mathcal{I}$ is $A$-identifiable. Hence, by Lemma~\ref{lemXYID}, if $|\mathcal{I}|\le 2$, then we require at most two vertices from $A$ to separate and dominate the vertices in $\mathcal{I}$. Moreover, we have $|A|\geq3$ due to the assumption $\deg(u)+\deg(v)\geq5$. Notice that set $A_{uv}\setminus\{z_1,z_2\}$ remains an identifying code of $G_{AI}$ when $z_1\in N_u$ and $z_2\in N_v$ and they are not required for dominating or separating vertices in $\mathcal{I}$. Indeed, these vertices will be the only ones which are adjacent to exactly $u$ or $v$ in the identifying code.

Hence, if $|\mathcal{I}|=2$ and  $|A_{uv}|=2\Delta-2-a\leq2\Delta-2$ for some $a\geq0$, then $A_{uv}$ is an identifying code of $G_{AI}$ that satisfies the conjectured bound for $G_{AI}$. Indeed,
\[
|A_{uv}|= 2\Delta-2-a= \left( \frac{\Delta-1}{\Delta} \right) ((2\Delta-2)+2)-a \le \left( \frac{\Delta-1}{\Delta} \right) n(G_{AI}).
\]
If $|\mathcal{I}|=2$ and $|A_{uv}|\geq2\Delta-1$, then we may remove one well-chosen vertex from $A$  since we require at most two vertices of $A$ to identify vertices in $\mathcal{I}$ while $|A|\geq3$, and the resulting set remains an identifying code (as described above) while satisfying the conjectured bound in $G_{AI}$. Indeed, $2\Delta-2< \left( \frac{\Delta - 1}{\Delta} \right)((2\Delta-2)+3)$ and $2\Delta-1< \left( \frac{\Delta - 1}{\Delta} \right)((2\Delta-1)+3)$.
If $|\mathcal{I}|=1$ and $|A_{uv}|\le 2\Delta-1$, then we may remove one well-chosen vertex from $N_u$. If $|\mathcal{I}|=1$ and $|A_{uv}|= 2\Delta$, then we may remove one well-chosen vertex from $N_u$ and $N_v$. Finally, if $\mathcal{I}=\emptyset$, then we may just remove any single vertex from both $N_u$ and $N_v$ to obtain an identifying code satisfying conjectured bound.

We are now ready to complete the proof of Claim~\ref{c:claimB}.
Observe that we have obtained above an identifying code, which contains $u$ and $v$, satisfying the conjectured bound in $G_{AI}$. Denote by $\mathcal{F}_{uv}$ the components in $G_{uv}$. Furthermore, observe that if $F\in \mathcal{F}_{uv}$ has size at least $2$ in $G_{uv}$, then $F$ is also a component in $G-G_{AI}$. Furthermore, by Claims~\ref{c:claimB.3}, \ref{c:claimB.4}, \ref{c:claimB.5}, \ref{c:claimB.6}, \ref{c:claimB.7} and \ref{c:claimB.8}, every component $F$ of $\mathcal{F}_{uv}$ admits an identifying code and we have $F\notin \mathcal{F}_\Delta$, and so, by the induction hypothesis, $\ID(F) \le \left( \frac{\Delta - 1}{\Delta} \right)n(F)$. Let us denote by $C_F$ an optimal identifying code in component $F$ and by $C_{AI}$ an identifying code in $G_{AI}$ which contains $u$ and $v$, and has size at most $\left( \frac{\Delta - 1}{\Delta}\right)n(G_{AI})$, as constructed above. Since there are no edges between components of $\mathcal{F}_{uv}$ and each edge from such a component $F$ is to a vertex in $A$, which are separated from vertices in $F$ by $u$ and $v$, we have an identifying code
\[
C=C_{AI}\cup\bigcup_{F\in \mathcal{F}_{uv}, n(F)\geq3} C_F
\]
and
\begin{align*}|C|&\le \left( \frac{\Delta - 1}{\Delta}\right)n(G_{AI})+\sum_{F\in \mathcal{F}_{uv}, n(F)\geq3}\left( \frac{\Delta - 1}{\Delta}\right)n(F)\\
&=\left( \frac{\Delta - 1}{\Delta}\right)n.\end{align*}

This completes the proof of Claim~\ref{c:claimB}.~\smallqed

\medskip

\begin{figure}[!htb]
\centering
\centering
\begin{tikzpicture}[
blacknode/.style={circle, draw=black!, fill=black!, thick},
whitenode/.style={circle, draw=black!, fill=white!, thick},
scale=1]
\node[blacknode,label={90:$u$}] (u) at (0,0) {};
\node[whitenode,label={90:$v$}] (v) at (3,0) {};

\node[whitenode] (2) at (-2,-2) {};
\node (2a) at (-1,-2) {};
\node (2b) at (0,-2) {};
\node[whitenode,label={160:$u'$}] (3) at (1,-2) {};
\node at ($(2)+(-0.75,0)$) {$N_u$};
\node at ($(2)+(1.5,0)$) {$\cdots$};
\draw[rounded corners] ($(2)+(-1.25,0.75)$) rectangle ($(3)+(0.5,-0.75)$);

\node[blacknode,label={20:$v'$}] (4) at (2.5,-2) {};
\node[blacknode] (4a) at (4,-2) {};
\node[whitenode] (4b) at (6.5,-2) {};
\node[whitenode] (5) at (8,-2) {};
\node at ($(5)+(0.75,0)$) {$N_v$};
\node at ($(4)+(0.75,0)$) {$\cdots$};
\node at ($(5)+(-0.75,0)$) {$\cdots$};
\draw[rounded corners] ($(4)+(-0.75,0.75)$) rectangle ($(5)+(1.25,-0.75)$);

\node at ($(4)+(0.75,-0.25)$) {$X$};
\draw[rounded corners] ($(4)+(-0.5,0.5)$) rectangle ($(4a)+(0.5,-0.5)$);

\draw[rounded corners] ($(2)+(-1.5,1)$) rectangle ($(5)+(2.25,-1)$);
\node at ($(5)+(1.75,0)$) {$A$};

\draw[rounded corners] ($(2)+(-1.75,2.75)$) rectangle ($(5)+(2.5,-1.25)$);
\node at ($(v)+(6.75,0.25)$) {$A_{uv}$};

\draw[rounded corners] ($(2)+(-1.75,-1.5)$) rectangle ($(5)+(2.5,-3)$);
\node at ($(5)+(1.75,-2.25)$) {$\barA_{uv}$};
\node at ($(5)+(3.25,-2.25)$) {$G_{uv}$};

\node[whitenode] (6a) at (1,-4.25) {};
\node[whitenode] (6b) at ($(6a)+(4,0)$) {};
\node at ($(6a)+(2,0)$) {$\cdots$};

\draw[-, thick,dashed] (u) -- (v) node[midway,above] {$e$};
\draw[-, thick] (u) -- (2);
\draw[-, thick] (4b) -- (v) -- (4a);
\draw[-, thick] (u) edge ($(u)!0.35!(2a)$) edge [dotted] ($(u)!0.45!(2a)$);
\draw[-, thick] (u) edge ($(u)!0.35!(2b)$) edge [dotted] ($(u)!0.45!(2b)$);
\draw[-, thick] (2) edge ($(2)!0.7!(6a)$) edge [dotted] ($(2)!0.8!(6a)$);
\draw[-, thick] (3) edge ($(3)!0.7!(6a)$) edge [dotted] ($(3)!0.8!(6a)$);
\draw[-, thick] (4a) edge ($(4a)!0.65!(6b)$) edge [dotted] ($(4a)!0.75!(6b)$);
\draw[-, thick] (4b) edge ($(4b)!0.65!(6a)$) edge [dotted] ($(4b)!0.75!(6a)$);
\draw[-, thick] (5) edge ($(5)!0.7!(6b)$) edge [dotted] ($(5)!0.8!(6b)$);
\draw[-, thick] (2) edge ($(2)!0.7!(6b)$) edge [dotted] ($(2)!0.8!(6b)$);
\draw[-, thick] (4) edge ($(4)!0.7!(6b)$) edge [dotted] ($(4)!0.8!(6b)$);

\draw[-, thick] (3) to[bend right=45] (4a);
\draw[-, thick] (u) -- (3) -- (4);
\draw[-, thick] (v) -- (4);
\draw[-, thick] (v) -- (5);
\end{tikzpicture}
\caption{The setting of the proof of Theorem~\ref{thm:trianglefree} after applying Claim~\ref{c:claimB}. Here, $\{u',v\}$ is not identified by $C'$ (black vertices).}
\label{fig:general_setting_after_ClaimB}
\end{figure}
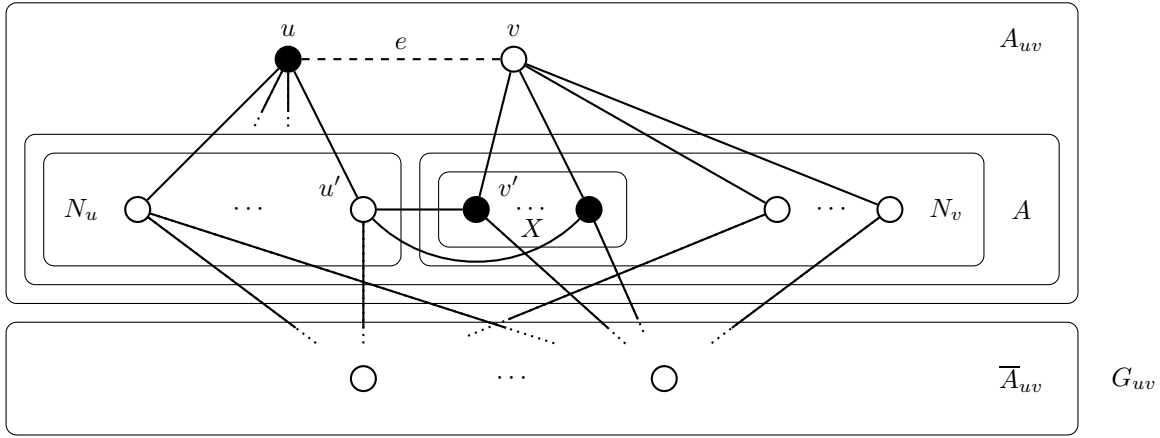

By Claim~\ref{c:claimB}, we may assume that there exists an optimal identifying code $C'$ in $G'$ such that $\{u',v\}$ or $\{u,v'\}$ is not identified by $C'$ in $G$, where $u'$ is some neighbor of $u$ different from $v$ and $v'$ is some neighbor of $v$ different from $u$.
Necessarily, in the case where $\{u',v\}$ is not identified, $u \in C'$ and $v \notin C'$, and if $X = N_{G'}(v) \cap C'$, then $X \ne \emptyset$ (in order for $C'$ to dominate~$v$ in $G'$) and $N_{G'}(u') \cap C' = X \cup \{u\}$. (We have the symmetric facts for the case where $\{u,v'\}$ is not identified.) Let $v'$ be an arbitrary vertex in $X$, and let $Q_{uv}$ be the cycle $uu'v'vu$ in $G$. Recall that $G_{uv} = G - N_G[u] - N_G[v]$. Further recall that $N_u = N_G(u) \setminus \{v\}$, $N_v = N_G(v) \setminus \{u\}$, and $A = N_u \cup N_v$. See Figure~\ref{fig:general_setting_after_ClaimB} for an illustration.

\begin{claim}
\label{c:claimC}
If $F$ is a component in $G_{uv}$ and $F \in \mathcal{F}_{\Delta}$, then
\[
\ID(G) \le \left( \frac{\Delta - 1}{\Delta}\right) n.
\]
\end{claim}
\proof Suppose that $F$ is a component in $G_{uv}$ and $F \in \mathcal{F}_{\Delta}$. Suppose, firstly, that $F = K_{1,\Delta}$ where $\Delta \ge 3$. Let $V(F) = \{x,x_1,x_2,\ldots,x_{\Delta}\}$ where $x$ is the central vertex of $F$ with leaf neighbors $x_1, x_2, \ldots, x_{\Delta}$. Since $G$ is connected, there is an edge $f$ joining a vertex $z \in A$ and a vertex in $V(F)$. Renaming vertices if necessary, we may assume that $f = zx_\Delta$. Let $G^* = G - ( V(F) \setminus \{x_\Delta\})$. Let $G^*$ have order~$n^*$, and so $n^* = n - \Delta$. Further, let $\Delta(G^*) = \Delta^*$. Since at least one of $u$ and $v$ has degree at least~$3$ in $G$, and since the degrees of $u$ and $v$ are the same in $G$ and in $G^*$, we note that $\Delta^* \ge 3$. Thus, $G^*$ is a connected triangle-free graph and $\Delta \ge \Delta^* \ge 3$. Moreover, $G^*$ contains the cycle $Q_{uv}$ and has order~$n^* \ge 6$. These properties of $G^*$ imply that $G^* \notin \mathcal{F}_{\Delta^*}$. Applying the inductive hypothesis to $G^*$, we have
\[
\ID(G^*) \le \left( \frac{\Delta^* - 1}{\Delta^*}\right) n^* \le \left( \frac{\Delta - 1}{\Delta}\right) n^*=\left( \frac{\Delta - 1}{\Delta}\right)n-(\Delta-1).
\]

Let $C^*$ be an optimal identifying code in $G^*$, and so  $|C^*| = \ID(G^*)$. The code $C^* \cup \{x,x_1,\dots,x_{\Delta-2}\}$ is an  identifying code in $G$ if $\Delta\ge 4$, or $\Delta=3$ and if $N_G(x_1)\cap C^*\neq \emptyset$. Similarly, if $N_G(x_2)\cap C^*\neq \emptyset$ and $\Delta=3$, then $C^* \cup \{x,x_2\}$ is an identifying code in $G$. Finally, if $\Delta=3$ and  $N_G(x_1)\cap C^*=N_G(x_2)\cap C^*=\emptyset$, then $C^* \cup \{x_1,x_2\}$ is an identifying code in $G$. Hence, we have
\[
\begin{array}{lcl}
\ID(G) & \le & |C^*| + \Delta - 1 \1 \\
& = & \ID(G^*) + \Delta - 1 \1 \\
& \le & \displaystyle{ \left(\left( \frac{\Delta - 1}{\Delta}\right)n-\Delta+1\right)+\Delta-1} \1 \\
& = & \displaystyle{  \left( \frac{\Delta - 1}{\Delta}\right) n, }
\end{array}
\]
which yields the desired upper bound in the case when $F = K_{1,\Delta}$ where $\Delta \ge 3$. Hence, we may assume that $\Delta = 3$ and $F \in \mathcal{F}_{\Delta} \setminus \{K_{1,3}\}$. Thus, $F \in \{P_4,C_4,C_7\} \cup (\mathcal{T}_{3} \setminus \{K_{1,3}\})$.

We next distinguish two cases.

\smallskip

\noindent\textit{Case 1: $F \ne P_4$.} In this case, since the graph $G$ is connected and $\Delta = 3$, one can check that there exists an induced path $P \colon x_1x_2x_3$ in $F$ such that $G'' = G - V(P)$ is a connected graph. Let $n'' = n(G'')$, and so $n'' = n - 3$. Let $\Delta(G'') = \Delta''$. We note that $G''$ is a connected triangle-free graph and $\Delta \ge \Delta'' \ge 3$. Thus, $\Delta''=3$.
Moreover, $G''$ contains the cycle $Q_{uv}$ and has order~$n'' \ge 6$. In particular, $G'' \notin \mathcal{F}_{3}$. Applying the inductive hypothesis to $G''$, we have
\[
\ID(G'') \le  \frac{2}{3} n'' =\frac{2}{3}n-2.
\]

Let $C''$ be an optimal identifying code in $G''$, and so, $|C''| = \ID(G'')$. If $N_G[x_1]\cap C''=N_G[x_3]\cap C''=\emptyset$, then we consider $C'' \cup \{x_1,x_3\}$ as a potential identifying code in $G$. It suffices to check that $x_1,x_2,x_3$ are now identified by the new code. Note that with this code, each of $x_1$ and $x_3$ is only dominated by itself, and no other vertex is dominated only by $x_1$ or only by $x_3$. If there exists some vertex $y$ not separated from $x_2$ by $C''\cup\{x_1,x_3\}$, then $y$ must be adjacent to both $x_1,x_3$. Since $G$ is triangle-free, $y$ is not adjacent to $x_2$. Thus, $y\notin C''$ (otherwise $y,x_2$ would be separated). Thus, $y$ is dominated by its third neighbor $z$, which is in $C''$, and hence, $z$ is adjacent to $x_2$. Since $\Delta''=\Delta=3$, vertex $y$ has no other neighbors in $G$. Thus, $y\notin A$, as in that case, $y$ should be adjacent to $u$ or $v$, but $z\notin\{u,v\}$ since $z$ is adjacent to $x_2$. Hence, $y\in F$ and $F$ is $C_4$. Thus, if $F\neq C_4$, we are done. If $F= C_4$, we show that $C''\cup\{x_1,x_2\}$ is an identifying code in $G$. Indeed, $z,x_1,x_2$ are code vertices inducing a $P_3$ and we may apply Lemma~\ref{LemP3} on them. Furthermore, $x_3$ is separated from all other neighbors of $x_2$ since they are in the identifying code. Furthermore, $y$ is adjacent to $z$ and $x_1$. If this is true also for some vertex $w$, then $w\in A$ and $z$ is adjacent to $w,y,x_2$ which is not possible since $\Delta=3$ and $z$ is adjacent to $u$ or $v$. Hence, $C''\cup\{x_1,x_2\}$ is an identifying code in $G$.

Assume next that for an index $j\in \{1,3\}$, say $j=1$, we have $N_G[x_1]\cap C''\neq\emptyset$. Consider the set $C'' \cup \{x_1,x_2\}$. If it is an identifying code of $G$, then we are done. Thus, assume it is not the case. Vertices $x_1$, $x_2$ and the neighbor of $x_1$ in $C''$ induce a $P_3$ and thus by Lemma~\ref{LemP3}, they are uniquely identified. Hence, $x_3$ is not separated from some other vertex $w$. Thus, $w$ is a neighbor of $x_2$, and $w\notin C''$. As $w$ is dominated by $C''$, say by $w'$, vertices $w$ and $x_3$ have $w'\in C''$ as a second common neighbor.
Notice that in this case, vertices $x_3,x_2,w,w'$ form a 4-cycle and $N(x_2)=\{x_1,x_3,w\}$. However, now set $C'' \cup \{x_2,x_3\}$ is an identifying code in $G$. Indeed, $x_2,x_3$ and $w'$ are separated from other vertices by Lemma~\ref{LemP3}. Assume first that $w'$ is adjacent to $x_1$. Then, since $G''$ is connected, $n''\geq 6$ and $\Delta=3$, vertex $w$ has a neighbor other than $w'$ and $x_2$, say $w''$. Since $w$ and $w'$ are separated by $C''$, vertex $w''$ belongs to $C''$. Since $x_3$ and $w$ were not separated by $C''\cup\{x_1,x_2\}$, $w''$ is adjacent to $x_3$. Again since $n''\geq 6$, there is no edge between $x_1$ and $w''$ (otherwise all considered vertices have degree 3 and $n''=3$). Hence, $x_1$ is identified by $C''\cup\{x_2,x_3\}$. On the other hand, if $w'$ is not adjacent to $x_1$, then  vertex $w'$ separates $x_1$ from the two other neighbors of $x_2$. In both cases, $C''\cup\{x_2,x_3\}$ is an identifying code of $G$. Therefore,
\[
\begin{array}{lcl}
\ID(G) & \le & |C''| + 2 \1 \\
& = & \ID(G'') + 2 \1 \\
& \le & \displaystyle{  \left(\frac{2}{3}n-2\right)  + 2 } \1 \\
& = & \displaystyle{   \frac{2}{3} n, }
\end{array}
\]
which yields the desired upper bound.

\smallskip

\noindent\textit{Case 2: $F = P_4$.} Let $F$ be the path $x_1x_2x_3x_4$. If $x_4$ is adjacent to a vertex in the set $A$, then as before, there exists an induced path $P \colon x_1x_2x_3$ in $F$ such that $G'' = G - V(P)$ is a connected graph of order at least 6 not in $\mathcal{F}_3$. Then, the same arguments as in Case~1 apply and we obtain an identifying code of the desired size.

Hence, we may assume that $x_4$ has degree~$1$ in $G$ (with $x_3$ as its unique neighbor in $G$). Analogously, we may assume that $x_1$ has degree~$1$ in $G$  (with $x_2$ as its unique neighbor in $G$). Since $G$ is connected, we may assume, renaming vertices if necessary, that $x_2$ is adjacent to a vertex in the set $A$.

We now consider the graph $G_F = G - \{x_3,x_4\}$. Let $n_F = n(G_F)$, and so $n_F = n - 2$. Let $\Delta(G_F) = \Delta_F$. We note that $G_F$ is a connected triangle-free graph and $\Delta \ge \Delta_F \ge 3$. Thus, $\Delta_F=3$. 
Moreover, $G_F$ contains the cycle $Q_{uv}$ and has order~$n_F \ge 7$. In particular, $G_F \notin \mathcal{F}_{3}$. Applying the inductive hypothesis to $G_F$, we have
\[
\ID(G_F) \le \frac{2}{3} n_F =\frac{2}{3}(n-2).
\]

Let $C_F$ be an optimal identifying code in $G_F$, and so $|C_F| = \ID(G_F)$. If $x_2 \in C_F$, then in order to identify the vertices $x_1$ and $x_2$, the code $C_F$ contains at least one neighbor of $x_2$ in $G_F$. In this case, we let $C = C_F \cup \{x_3\}$. If $x_2 \notin C_F$, then $x_1 \in C_F$ and at least one neighbor of $x_2$ in the set $A$ belongs to $C_F$. In this case, we let $C = (C_F \setminus \{x_1\}) \cup \{x_2,x_3\}$. In both cases (using Lemma~\ref{LemP3} in the second case), $C$ is an identifying code of $G$ and $|C| = |C_F| + 1$, implying that

\[
\begin{array}{lcl}
\ID(G) \le |C| & = & \ID(G_F) + 1 \1 \\
& \leq & \displaystyle{  \frac{2}{3}(n-2)+1 } \1 \\
& < & \displaystyle{   \frac{2}{3} n. }
\end{array}
\]
This completes the proof of Claim~\ref{c:claimC}.~\smallqed

\medskip

By Claim~\ref{c:claimC}, we will from now on assume that if $F$ is a component in $G_{uv}$ of order at least~$3$, then $F \notin \mathcal{F}_{\Delta}$.

\medskip
Let $B$ be the set of all vertices that belong to a $P_1$-component or to a $P_2$-component in $G_{uv}$. Furthermore, let us denote $G^*=G[A\cup B\cup\{u,v\}]$. Note that it is possible that $G^*=G$ and hence we do not apply the induction hypothesis directly to $G^*$. Recall that $A = N_u \cup N_v$ and let $|A| = a$. Further, let $|B|= b$. Thus, $V(G^*) = A \cup B \cup \{u,v\}$ and $n^* = a + b + 2$. Next, in Claim~\ref{claimGStar}, we will show that $G^*$ admits an identifying code of cardinality at most $\left( \frac{\Delta-1}{\Delta} \right) n^*$ that contains both vertices $u$ and $v$. Notice that $G^*\not\in \mathcal{F}_\Delta$ since $G^*$ contains a 4-cycle and a vertex of degree at least~$3$.

 \begin{claim}\label{claimGStar}
 Graph $G^*$ admits an identifying code containing vertices $u$ and $v$ of cardinality at most $$\left(\frac{\Delta-1}{\Delta}\right) n^*.$$
 \end{claim}
 \proof
 Since $G$ is a connected graph, every vertex in a $P_1$-component of $G[B]$ is adjacent in $G^*$ to at least one vertex in~$A$. Moreover, by the connectivity of $G^*$, at least one vertex from every $P_2$-component of $G[B]$ is adjacent in $G^*$ to at least one vertex in~$A$. However, it is possible that a $P_2$-component of $G[B]$ contains exactly one vertex that has degree~$1$ in $G^*$ (and is therefore not adjacent in $G^*$ to a vertex from the set $A$). Let $B_A$ be the set of all vertices in $B$ that have at least one neighbor in $A$ in the graph $G^*$, and let $B_L = B \setminus B_A$. Thus, each vertex in $B_L$ is a vertex of degree~$1$ in $G^*$ that belongs to a  $P_2$-component of $G[B]$ and is adjacent to no vertex of $A$ in $G^*$. Possibly, $B_L = \emptyset$.

Let $B = (B_1,\ldots,B_k)$ be a partition of the set $B$ such that the following properties hold. \\ [-18pt]
\begin{enumerate}
\item[{\rm (a)}] All vertices in the set $B_i \cap B_A$ have the same neighborhood in the set $A$ for all $i \in [k]$, that is, if $x,y \in B_i \cap B_A$, then $N_{G^*}(x) \cap A = N_{G^*}(y) \cap A \ne \emptyset$. \1
\item[{\rm (b)}] Each vertex in the set $B_i \cap B_L$ has its unique neighbor (in ${G^*}$) in the set $B_i \cap B_A$ for all $i \in [k]$.\1
\item[{\rm (c)}] Vertices in distinct sets $B_i \cap B_A$ and $B_j \cap B_A$ have different neighborhoods in $A$, that is, if $x \in B_i \cap B_A$ and $y \in B_j \cap B_A$ where $1 \le i < j \le k$, then $N_{G^*}(x) \cap A \ne N_{G^*}(y) \cap A$.
\end{enumerate}

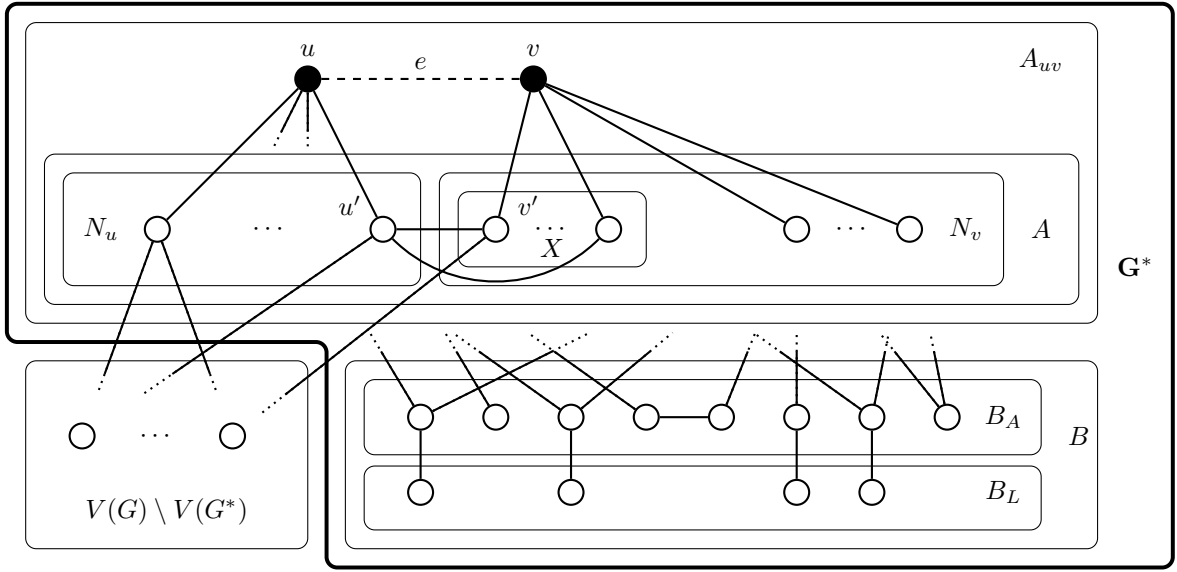
\begin{figure}[!htb]
\centering
\centering
\begin{tikzpicture}[
blacknode/.style={circle, draw=black!, fill=black!, thick},
whitenode/.style={circle, draw=black!, fill=white!, thick},
scale=1]

\node[blacknode,label={90:$u$}] (u) at (0,0) {};
\node[blacknode,label={90:$v$}] (v) at (3,0) {};

\node[whitenode] (2) at (-2,-2) {};
\node (2a) at (-1,-2) {};
\node (2b) at (0,-2) {};
\node[whitenode,label={160:$u'$}] (3) at (1,-2) {};
\node at ($(2)+(-0.75,0)$) {$N_u$};
\node at ($(2)+(1.5,0)$) {$\cdots$};
\draw[rounded corners] ($(2)+(-1.25,0.75)$) rectangle ($(3)+(0.5,-0.75)$);

\node[whitenode,label={20:$v'$}] (4) at (2.5,-2) {};
\node[whitenode] (4a) at (4,-2) {};
\node[whitenode] (4b) at (6.5,-2) {};
\node[whitenode] (5) at (8,-2) {};
\node at ($(5)+(0.75,0)$) {$N_v$};
\node at ($(4)+(0.75,0)$) {$\cdots$};
\node at ($(5)+(-0.75,0)$) {$\cdots$};
\draw[rounded corners] ($(4)+(-0.75,0.75)$) rectangle ($(5)+(1.25,-0.75)$);

\node at ($(4)+(0.75,-0.25)$) {$X$};
\draw[rounded corners] ($(4)+(-0.5,0.5)$) rectangle ($(4a)+(0.5,-0.5)$);

\draw[rounded corners] ($(2)+(-1.5,1)$) rectangle ($(5)+(2.25,-1)$);
\node at ($(5)+(1.75,0)$) {$A$};

\draw[rounded corners] ($(2)+(-1.75,2.75)$) rectangle ($(5)+(2.5,-1.25)$);
\node at ($(v)+(6.75,0.25)$) {$A_{uv}$};

\draw[rounded corners] ($(2)+(-1.75,-1.75)$) rectangle ($(2)+(2,-4.25)$);
\node at ($(2)+(0.125,-3.75)$) {$V(G)\setminus V(G^*)$};
\node[whitenode] (6a) at ($(2)+(-1,-2.75)$) {};
\node[whitenode] (6b) at ($(6a)+(2,0)$) {};
\node at ($(6a)+(1,0)$) {$\cdots$};

\draw[rounded corners] ($(2)+(2.5,-1.75)$) rectangle ($(5)+(2.5,-4.25)$);
\node at ($(5)+(2.25,-2.75)$) {$B$};
\draw[rounded corners] ($(2)+(2.75,-2)$) rectangle ($(5)+(1.75,-3)$);
\node at ($(5)+(1.25,-2.5)$) {$B_A$};
\draw[rounded corners] ($(2)+(2.75,-3.15)$) rectangle ($(5)+(1.75,-4)$);
\node at ($(5)+(1.25,-3.5)$) {$B_L$};
\node[whitenode] (b1) at ($(2)+(3.5,-2.5)$) {};
\node[whitenode] (b1p) at ($(b1)+(0,-1)$) {};
\node[whitenode] (b2) at ($(b1)+(1,0)$) {};
\node[whitenode] (b3) at ($(b2)+(1,0)$) {};
\node[whitenode] (b3p) at ($(b3)+(0,-1)$) {};
\node[whitenode] (b4) at ($(b3)+(1,0)$) {};
\node[whitenode] (b5) at ($(b4)+(1,0)$) {};
\node[whitenode] (b6) at ($(b5)+(1,0)$) {};
\node[whitenode] (b6p) at ($(b6)+(0,-1)$) {};
\node[whitenode] (b7) at ($(b6)+(1,0)$) {};
\node[whitenode] (b7p) at ($(b7)+(0,-1)$) {};
\node[whitenode] (b8) at ($(b7)+(1,0)$) {};
\draw[-, thick] (b1) edge ($(b1)!0.35!(2b)$) edge [dotted] ($(b1)!0.45!(2b)$);
\draw[-, thick] (b1) edge ($(b1)!0.35!(4b)$) edge [dotted] ($(b1)!0.45!(4b)$);
\draw[-, thick] (b2) edge ($(b2)!0.35!(3)$) edge [dotted] ($(b2)!0.45!(3)$);
\draw[-, thick] (b3) edge ($(b3)!0.35!(2b)$) edge [dotted] ($(b3)!0.45!(2b)$);
\draw[-, thick] (b3) edge ($(b3)!0.35!(4b)$) edge [dotted] ($(b3)!0.45!(4b)$);
\draw[-, thick] (b5) edge ($(b5)!0.35!(4b)$) edge [dotted] ($(b5)!0.45!(4b)$);
\draw[-, thick] (b4) edge ($(b4)!0.35!(3)$) edge [dotted] ($(b4)!0.45!(3)$);
\draw[-, thick] (b6) edge ($(b6)!0.35!(4b)$) edge [dotted] ($(b6)!0.45!(4b)$);
\draw[-, thick] (b7) edge ($(b7)!0.35!(4a)$) edge [dotted] ($(b7)!0.45!(4a)$);
\draw[-, thick] (b7) edge ($(b7)!0.35!(5)$) edge [dotted] ($(b7)!0.45!(5)$);
\draw[-, thick] (b8) edge ($(b8)!0.35!(4b)$) edge [dotted] ($(b8)!0.45!(4b)$);
\draw[-, thick] (b8) edge ($(b8)!0.35!(5)$) edge [dotted] ($(b8)!0.45!(5)$);

\node at ($(5)+(3,-0.5)$) {$\mathbf{G^*}$};
\draw[rounded corners,line width=1.5pt] ($(u)+(0,1)$) -- ++(11.5,0) -- ++(0,-7.5) -- ++(-11.25,0) -- ++(0,3) -- ++(-4.25,0) -- ++(0,4.5) -- ($(u)+(0,1)$);

\draw[-, thick,dashed] (u) -- (v) node[midway,above] {$e$};
\draw[-, thick] (u) -- (2);
\draw[-, thick] (4b) -- (v) -- (4a);
\draw[-, thick] (u) edge ($(u)!0.35!(2a)$) edge [dotted] ($(u)!0.45!(2a)$);
\draw[-, thick] (u) edge ($(u)!0.35!(2b)$) edge [dotted] ($(u)!0.45!(2b)$);
\draw[-, thick] (2) edge ($(2)!0.7!(6a)$) edge [dotted] ($(2)!0.8!(6a)$);
\draw[-, thick] (3) edge ($(3)!0.7!(6a)$) edge [dotted] ($(3)!0.8!(6a)$);
\draw[-, thick] (2) edge ($(2)!0.7!(6b)$) edge [dotted] ($(2)!0.8!(6b)$);
\draw[-, thick] (4) edge ($(4)!0.8!(6b)$) edge [dotted] ($(4)!0.9!(6b)$);

\draw[-, thick] (3) to[bend right=45] (4a);
\draw[-, thick] (u) -- (3) -- (4);
\draw[-, thick] (v) -- (4);
\draw[-, thick] (v) -- (5);
\draw[-, thick] (b1) -- (b1p) (b3) -- (b3p) (b6) -- (b6p) (b7) -- (b7p) (b4)--(b5);

\end{tikzpicture}
\caption{The setting of Claim~\ref{claimGStar} in the proof of Theorem~\ref{thm:trianglefree}. The aim is to construct a small identifying code containing both $u$ and $v$.}
\label{fig:general_setting_clmD}
\end{figure}

An illustration is given in Figure~\ref{fig:general_setting_clmD}. We note that $|B_i| = |B_i \cap B_A| + |B_i \cap B_L|$ for all $i \in [k]$. Moreover, $|B_i \cap B_L| \le |B_i \cap B_A|$ and $|B_i \cap B_A| \ge 1$ for all $i \in [k]$. If $|B_i \cap B_L| = |B_i \cap B_A|$, then we observe that the subgraph $G[B_i]$ induced by the set $B_i$ consists of vertex-disjoint copies of $K_2$ (and each such copy of $K_2$ contains exactly one vertex of degree~$1$ in $G$). Since each vertex in $A$ is adjacent to either the vertex~$u$ or~$v$, each vertex in $A$ has at most~$\Delta - 1$ neighbors in the set $B$, implying that $|B_i \cap B_A| \le \Delta - 1$ for all $i \in [k]$. Let
\[
|B_i| = b_i \1
\]
for all $i \in [k]$. Thus, if $B_i \cap B_L = \emptyset$, then $b_i = |B_i \cap B_A| \le \Delta - 1$, while if $B_i \cap B_L \ne \emptyset$, then $b_i \le 2|B_i \cap B_A| \le 2(\Delta - 1)$ for all $i \in [k]$. Recall that
\[
b = |B| = \sum_{i=1}^k b_i.
\]

Since $G$ is triangle-free, we note that $B_i \cap B_A$ is an independent set in ${G^*}$ for all $i \in [k]$. Thus, if two vertices $w$ and $z$ belong to the same $P_2$-component in $G[B]$, then either $w \in B_i$ and $z \in B_j$ where $1 \le i, j \le k$ and $i \ne j$, or one of $w$ and $z$ belongs to $B_i \cap B_A$ and the other to $B_i \cap B_L$ for some $i \in [k]$. For $i \in [k]$, we now define the set $B_i^*$ as follows. If $|B_i \cap B_L| = |B_i \cap B_A|$, then we define
\[
B_i^* = B_i \cap B_L,
\]
while if $|B_i \cap B_L| < |B_i \cap B_A|$, then we let $w_i$ be an arbitrary vertex in the set $B_i \cap B_A$ that is isolated in the subgraph $G[B_i]$ induced by $B_i$ (and so, $w_i$ has no neighbor in $B_i$), and we define
\[
B_i^* = (B_i \cap B_L) \cup \{w_i\}.
\]

We now define
\[
B^* = \bigcup_{i=1}^k B_i^*. 
\]

Moreover, let $|B_i^*| = b_i^*$ for all $i \in [k]$ and
\[
b^*=\sum_{i=1}^k b_i^*.
\]

\begin{figure}[!htb]
\centering
\centering
\begin{tikzpicture}[
blacknode/.style={circle, draw=black!, fill=black!, thick},
whitenode/.style={circle, draw=black!, fill=white!, thick},
scale=1]

\node[whitenode] (b1) at (0,0) {};
\node[whitenode,fill=black!20] (b1p) at ($(b1)+(0,-1)$) {};
\node[whitenode] (b2) at ($(b1)+(0.5,0)$) {};
\node[whitenode,fill=black!20] (b2p) at ($(b2)+(0,-1)$) {};
\node[whitenode] (b3) at ($(b2)+(0.5,0)$) {};
\node[whitenode,fill=black!20] (b3p) at ($(b3)+(0,-1)$) {};
\draw[rounded corners] ($(b1)+(-0.25,0.25)$) rectangle ($(b3p)+(0.25,-0.25)$);
\node at ($(b2p)+(0,-0.75)$) {$B_1$};

\node[whitenode] (b4) at ($(b3)+(1,0)$) {};
\node[whitenode,fill=black!20] (b4p) at ($(b4)+(0,-1)$) {};
\node[whitenode,fill=black!20] (b5) at ($(b4)+(0.5,0)$) {};
\node[whitenode] (b6) at ($(b5)+(0.5,0)$) {};
\draw[rounded corners] ($(b4)+(-0.25,0.25)$) rectangle ($(b6)+(0.25,-1.25)$);
\node at ($(b5)+(0,-1.75)$) {$B_2$};

\node[whitenode,fill=black!20] (b7) at ($(b6)+(1,0)$) {};
\node[whitenode] (b8) at ($(b7)+(0.5,0)$) {};
\node[whitenode] (b9) at ($(b8)+(0.5,0)$) {};
\node[whitenode,fill=black!20] (b9p) at ($(b9)+(0,-1)$) {};
\node[whitenode] (b10) at ($(b9)+(0.5,0)$) {};
\draw[rounded corners] ($(b7)+(-0.25,0.25)$) rectangle ($(b10)+(0.25,-1.25)$);
\node at ($(b8)+(0.25,-1.75)$) {$B_3$};

\node[whitenode] (b11) at ($(b10)+(1,0)$) {};
\node[whitenode,fill=black!20] (b11p) at ($(b11)+(0,-1)$) {};
\node[whitenode] (b12) at ($(b11)+(0.5,0)$) {};
\node[whitenode,fill=black!20] (b12p) at ($(b12)+(0,-1)$) {};
\node[whitenode,fill=black!20] (b13) at ($(b12)+(0.5,0)$) {};
\draw[rounded corners] ($(b11)+(-0.25,0.25)$) rectangle ($(b13)+(0.25,-1.25)$);
\node at ($(b12)+(0,-1.75)$) {$B_4$};

\node[whitenode,fill=black!20] (b14) at ($(b13)+(1,0)$) {};
\node[whitenode] (b15) at ($(b14)+(0.5,0)$) {};
\draw[rounded corners] ($(b14)+(-0.25,0.25)$) rectangle ($(b15)+(0.25,-1.25)$);
\node at ($(b14)+(0.25,-1.75)$) {$B_5$};

\draw[-, thick] (b1)--(b1p) (b2)--(b2p) (b3)--(b3p) (b4)--(b4p) (b6)--(b7) (b9)--(b9p) (b11)--(b11p) (b12)--(b12p) (b13)--(b14);

\end{tikzpicture}
\caption{The construction of the set $B^*$ in the proof of Claim~\ref{claimGStar}. The vertices in $B^*$ are shaded.}
\label{fig:Bstar}
\end{figure}
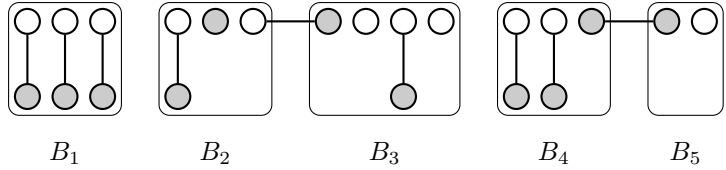

An illustration of the set $B^*$ is given in Figure~\ref{fig:Bstar}.

\medskip
\begin{subclaim}
\label{c:claimE}
The following hold. \\[-16pt]
\begin{enumerate}
\item[{\rm (a)}] $b_i^* \ge \frac{b_i}{\Delta - 1}$ for all $i \in [k]$. \1
\item[{\rm (b)}] $b \le (\Delta - 1)b^*$.
\end{enumerate}
\end{subclaim}
\proof By our earlier observations, if $B_i \cap B_L = \emptyset$, then $b_i \le \Delta - 1$ where $i \in [k]$. Moreover, $b_i \le 2(\Delta - 1)$ for all $i \in [k]$. If $b_i \le \Delta - 1$ for some $i \in [k]$, then
\[
\frac{b_i - b_i^*}{b_i} \le \frac{b_i - 1}{b_i} \le  \frac{\Delta - 2}{\Delta - 1}.
\]

We note that if $b_i \ge \Delta$ for some $i \in [k]$, then $b_i^* \ge 2$, and so in this case
\[
\frac{b_i - b_i^*}{b_i} \le \frac{b_i - 2}{b_i} \le \frac{2(\Delta - 1) - 2 }{2(\Delta - 1)} = \frac{\Delta - 2}{\Delta - 1}.
\]
Thus in both cases,
\[
\frac{b_i - b_i^*}{b_i} \le  \frac{\Delta - 2}{\Delta - 1}. \1
\]
Rearranging terms in the above inequality yields the inequality
\[
b_i^* \ge \frac{b_i}{\Delta - 1}
\]
for all $i \in [k]$. This proves property~(a) in the statement of the claim. Hence,
\[
b = \sum_{i=1}^k b_i \le (\Delta - 1) \sum_{i=1}^k b_i^* = (\Delta - 1)b^*,
\]
and so property~(b) in the statement of the claim holds.~\smallqed

\medskip

Let $z_i$ be an arbitrary vertex in $B_i \cap B_A$, and so $z_i \in B_i$ and $z_i$ has a neighbor in the set $A$ for all $i \in [k]$. Let $Z = \{z_1,z_2,\ldots,z_k\}$. By construction, the set $A$ identifies the set $Z$ since every pair of vertices in $Z$ have distinct neighborhoods in $A$.  Equivalently, the set $Z$ is $(Z,A)$-identifiable by the set $A$. Let $A^*$ be a subset of $A$ of minimum cardinality that $A$-identifies~$Z$. By Lemma~\ref{lemXYID}, $1 \le |A^*| \le k$. Since $b_i^* \ge 1$ for all $i \in [k]$, we note that $b^* \ge k$, and so $|A^*| \le b^*$. Recall that $Q_{uv} \colon uu'v'vu$ is a 4-cycle in $G$ and hence in $G^*$, where $u'$ is a neighbor of $u$ in $G$ different from $v$, and $v'$ is a neighbor of $v$ in $G$ different from $u$.

\begin{subclaim}
\label{c:claimF}
If $A = A^*$, then there exists an identifying code in $G^*$ containing vertices $u$ and $v$ with cardinality at most
\[
\left( \frac{\Delta - 1}{\Delta}\right) n^*.
\]
\end{subclaim}
\proof Suppose that $A = A^*$. By the minimality of the set $A^*$, every vertex in $A^*$ has a neighbor in $B$. Recall that $V(G^*) = A \cup B \cup \{u,v\}$ and $n^* = a + b + 2$. We now let
\[
C^* = A \cup \{u,v\} \cup (B \setminus B^*).
\] 
The set $C^*$ is an identifying code of $G^*$. Indeed, $C^*$ is a dominating set which is connected in $G^*$. Thus, by Lemma~\ref{LemP3}, it separates vertices in $C^*$. Furthermore, each vertex in $B^*$ has a unique neighborhood in $A\cup (B\setminus B^*)$.

\medskip
\begin{subsubclaim}
\label{c:claimF.1}
If $b^* \ge k + 2$, then  $|C^*| \le \left( \frac{\Delta - 1}{\Delta}\right) n^*.$
\end{subsubclaim}
\proof Suppose that $b^* \ge k + 2$. Thus, $a = |A^*| \le k \le b^* - 2$. By Claim~\ref{c:claimE}(b) and our earlier observations, we have
\[
\begin{array}{lcl}
 |C^*| & = & a + 2 + b - b^* \2 \\
& = & \displaystyle{ \left( \frac{\Delta - 1}{\Delta} \right) (a + 2 + b) + \frac{1}{\Delta}(a + 2 + b - \Delta b^*) } \2 \\
& \le & \displaystyle{ \left( \frac{\Delta - 1}{\Delta} \right)n^* + \frac{1}{\Delta}\big((b^* - 2) + 2 + (\Delta - 1)b^* - \Delta b^* \big) } \2 \\
& = & \displaystyle{ \left( \frac{\Delta - 1}{\Delta} \right)n^* }, \1 \\
\end{array}
\]
yielding the desired upper bound.~\smallqed

\medskip

By Claim~\ref{c:claimF.1}, we may assume that $b^* \le k + 1$, implying that $b^* \in \{k,k+1\}$. With this assumption, $b^*_i = 1$ for all $i \in [k]$, except for possibly exactly one value of~$i$ satisfying $b^*_i = 2$.

Recall that every vertex in $A^*$ has a neighbor in $B$. Thus, since $A = A^*$, we note in particular that both vertices $u'$ and $v'$ have a neighbor in $B$. Since $G^*$ is triangle-free, $u'$ and $v'$ do not have a common neighbor. Renaming the partitions in $B = (B_1,\ldots,B_k)$ if necessary, we may assume that $u'$ is adjacent to a vertex in $B_1$ and $v'$ is adjacent to a vertex in $B_2$. Since $u'$ is adjacent to both $u$ and $v'$, it has at most $\Delta - 2$ additional neighbors, implying that $|B_1 \cap B_A| \le \Delta - 2$. Analogously, $|B_2 \cap B_A| \le \Delta - 2$. Therefore,
\[
1 \le b_1^* \le b_1 \le 2(\Delta - 2) \hspace*{0.5cm} \mbox{and} \hspace*{0.5cm} 1 \le b_2^* \le b_2 \le 2(\Delta - 2).
\]
Let
\[
B_{3,k} = B \setminus (B_1 \cup B_2) \hspace*{0.5cm} \mbox{and} \hspace*{0.5cm} B_{3,k}^* = B^* \setminus (B_1^* \cup B_2^*). \1
\]

Further, let $b_{3,k} = |B_{3,k}|$ and $b_{3,k}^* = |B_{3,k}^*|$, and so $b = b_1 + b_2 + b_{3,k}$ and $b^* = b_1^* + b_2^* + b_{3,k}^*$. We note that
\[
b_{3,k} = \sum_{i=3}^k b_i \hspace*{0.5cm} \mbox{and} \hspace*{0.5cm} b_{3,k}^* = \sum_{i=3}^k b_i^* \ge \sum_{i=3}^k 1 = k-2.
\]

We now define a partition $(V_1,V_2,V_3)$ of $V(G^*)$ as follows. We let
\[
\begin{array}{lcl}
V_1 & = & B_1 \cup \{u,u'\}, \1 \\
V_2 & = & B_2 \cup \{v,v'\}, \1 \\
V_3 & = & B_{3,k} \cup (A \setminus \{u',v'\}). \1 \\
\end{array}
\]

Further we define $n_i = |V_i|$ for $i \in [3]$, and so $n^* = n_1 + n_2 + n_3$. We note that $n_1 = b_1 + 2$, $n_2 = b_2 + 2$, and $n_3 =  b_{3,k} + a - 2$. We again consider the identifying code $C^*$ of $G^*$
and we let $C_i^* = C^* \cap V_i$ for $i \in [3]$. Thus,
\[
|C^*| = \sum_{i=1}^3 |C_i^*|.
\]

\medskip
\begin{subsubclaim}
\label{c:claimF.2}
$|C_3^*| \le \left( \frac{\Delta - 1}{\Delta}\right) n_3$.
\end{subsubclaim}
\proof By Claim~\ref{c:claimE}(a), we have
\[
b_{3,k}^* = \sum_{i=3}^k b_i^* \ge \sum_{i=3}^k \frac{b_i}{\Delta - 1} = \frac{b_{3,k}}{\Delta - 1}, \1
\]
and so $(\Delta - 1)b_{3,k}^* \ge b_{3,k}$, or, equivalently,
\begin{equation}
\label{Eq:ClaimFa}
b_{3,k} - \Delta  b_{3,k}^* \le - b_{3,k}^* \le - k + 2.
\end{equation}

Recall that $n_3 =  a - 2 + b_{3,k}$ and that $a \le k$. Hence by Inequality~(\ref{Eq:ClaimFa}), we have
\[
\begin{array}{lcl}
|C_3^*| & = & |A \setminus \{u',v'\}| + b_{3,k} - b_{3,k}^* \2 \\
& = & a - 2 + b_{3,k} - b_{3,k}^* \2 \\
& = & \displaystyle{ \left( \frac{\Delta - 1}{\Delta} \right) (a - 2 + b_{3,k}) + \frac{1}{\Delta}(a - 2 + b_{3,k} - \Delta b_{3,k}^* ) } \2 \\
& \le & \displaystyle{ \left( \frac{\Delta - 1}{\Delta} \right) n_3 + \frac{1}{\Delta}(k - 2 - k + 2 ) } \2 \\
& = & \displaystyle{ \left( \frac{\Delta - 1}{\Delta} \right)n_3 }, \1 \\
\end{array}
\]
yielding the desired upper bound.~\smallqed

\medskip

\begin{subsubclaim}
\label{c:claimF.3}
If $b_i^* \ge \frac{b_i}{\Delta - 2}$, then $|C_i^*| \le \left( \frac{\Delta - 1}{\Delta}\right) n_i$ for $i \in [2]$.
\end{subsubclaim}
\proof Let $i \in [2]$. Recall that $b_i^* \ge 1$ and $n_i = b_i + 2$. Suppose that $b_i^* \ge \frac{b_i}{\Delta - 2}$. Thus, $(\Delta - 2)b_i^* \ge b_i$, or, equivalently, $b_i - \Delta b_i^* \le - 2b_i^* \le - 2$ noting that $b_i^* \ge 1$. By definition, we have $C_1^* = \{u,u'\} \cup (B_1 \setminus B_1^*)$ and $C_2^* = \{v,v'\} \cup (B_2 \setminus B_2^*)$. Hence,
\[
\begin{array}{lcl}
|C_i^*| & = & 2 + b_i - b_i^* \2 \\
& = & \displaystyle{ \left( \frac{\Delta - 1}{\Delta} \right) (2 + b_i) + \frac{1}{\Delta}(2 + b_i - \Delta b_i^* ) } \2 \\
& \le & \displaystyle{ \left( \frac{\Delta - 1}{\Delta} \right) n_i + \frac{1}{\Delta}(2 - 2 ) } \2 \\
& = & \displaystyle{ \left( \frac{\Delta - 1}{\Delta} \right)n_i }, \1 \\
\end{array}
\]
yielding the desired upper bound.~\smallqed

\medskip

If $b_i^* \ge \frac{b_i}{\Delta - 2}$ for $i \in [2]$, then by Claims~\ref{c:claimF.2} and~\ref{c:claimF.3}, we have
\[
 |C^*| = \sum_{i=1}^3 |C_i^*| \le \sum_{i=1}^3 \left( \frac{\Delta - 1}{\Delta}\right) n_i = \left( \frac{\Delta - 1}{\Delta}\right) n^*, \2
\]
yielding the desired upper bound. Hence, we may assume that $b_i^* < \frac{b_i}{\Delta - 2}$ for some $i \in [2]$. By symmetry, we may assume that $b_1^* < \frac{b_1}{\Delta - 2}$. Let
\[
|B_1 \cap B_A| = t_1 + t_2 \hspace*{0.5cm} \mbox{and} \hspace*{0.5cm}  |B_1 \cap B_L| = t_2.
\]

We note that $1 \le t_1 + t_2 \le \Delta - 2$.

If $t_1 \ge 1$ and $t_2 \ge 1$, then $b_1 = t_1 + 2t_2$ and $b_1^* = 1 + t_2 \ge 2$, implying that
\[
\frac{b_1}{\Delta - 2} = \frac{t_2 + (t_1 + t_2)}{\Delta - 2}  \le \frac{t_2 + (\Delta - 2)}{\Delta - 2} = \frac{t_2}{\Delta - 2} + 1 \le \frac{\Delta - 3}{\Delta - 2} + 1 < 2 \le b_1^*. \1
\]

If $t_1 \ge 1$ and $t_2  = 0$, then $b_1 = t_1 \le \Delta - 2$ and $b_1^* = 1$, implying that
\[
\frac{b_1}{\Delta - 2} \le \frac{\Delta - 2}{\Delta - 2} = 1 = b_1^*. \1
\]

If $t_1 = 0$ and $t_2 \ge 2$, then $b_1 = 2t_2$ and $b_1^* = t_2 \ge 2$. Moreover, $t_2 \le \Delta - 2$. Thus,
\[
\frac{b_1}{\Delta - 2} = \frac{2t_2}{\Delta - 2}  \le \frac{2(\Delta - 2)}{\Delta - 2} = 2 \le b_1^*. \1
\]

If $t_1 = 0$, $t_2 = 1$ and $\Delta \ge 4$, then $b_1 = 2$, $b_1^* = 1$ and
\[
\frac{b_1}{\Delta - 2} = \frac{2}{\Delta - 2}  \le \frac{2}{2} = 1 = b_1^*. \1
\]

In all the above four cases, we contradict our assumption that $b_1^* < \frac{b_1}{\Delta - 2}$. Hence, $t_1 = 0$, $t_2 = 1$, and $\Delta = 3$. In this case, $B_1$ induces a $P_2$-component in $G^*$. Let $B_1 \cap B_A = \{u_1\}$ and let $B_1 \cap B_L = \{u_2\}$, and so $uu'u_1u_2$ is a path in $G^*$. We note that $u_2$ is a vertex of degree~$1$ in $G^*$. Since $\Delta = 3$ and $v'$ is adjacent to $v$ and $u'$, we note that either $B_2 = P_1$ or $B_2 = P_2$. We denote the vertex in $B_2\cap N(v')$ by $v_1$ and if $v_1$ has another adjacent vertex outside of $A$, we denote it by $v_2$. By our choice of $u$ and $v$ (maximizing $\deg_G(u) + \deg_G(v)$ for all adjacent $u,v$ with $uv$ a cycle edge), we note that both $u$ and $v$ have degree~$3$ in $G^*$ (otherwise, $\deg_G(u') + \deg_G(v') > \deg_G(u)+ \deg_G(v)$). We denote the third neighbor of $u$ and $v$ by $u''$ and $v''$, respectively.

Let us assume next that $|N(u_1)\cap A|=2$. Since $G$ is triangle-free, the other vertex in $N(u_1)\cap A$ is not~$v'$. Notice that in this case, we may remove $u'$ from $C_1^*$, resulting in $C_1^{**}$, and the set $C^{**}=C_1^{**}\cup C_2^{*}\cup C_3^{*}$ is an identifying code in $G^*$. Indeed, $C^*$ was an identifying code in $G^*$. Moreover, every neighbor of $u'$ is in $C^{**}$ and vertices in $C^{**}$ induce a single component in $G[C^{**}]$. Hence, every vertex in $C^{**}$ is separated by Lemma~\ref{LemP3} while $u'$ is the unique vertex not in $C^{**}$ adjacent to vertices $u$ and $u_1$. Notice that $|C^{**}|\le \frac{2}{3}n^*$. Indeed, we have $n_1=4$,  $|C_1^{**}|=2$, $(n_2,|C_2^{*}|)\in\{(3,2),(4,3)\}$ and $|C_3^*|\le \frac{2}{3}n_3$ by Claim~\ref{c:claimF.2}. Since $\frac{|C_1^{**}|+|C_2^*|}{n_1+n_2}\le \frac{5}{8}<\frac{2}{3}$, we have $|C^{**}|\le \frac{2}{3}n^*$. Therefore, we may assume that one of $B_1$ or $B_2$ is a $P_2$-component such that the vertex in $B_A$ is adjacent to exactly one vertex in $A$. We may assume from now on without loss of generality, that $N(u_1)\cap A=\{u'\}$.

Recall that in graph $G'$, any minimum-size identifying code $C'$ is such that it does not separate either vertices $\{u,v'\}$ or $\{v,u'\}$ in $G$. Thus, $C'\cap\{u,v,v',u'\}$ is either $\{u,v'\}$ or $\{v,u'\}$. Furthermore, since $u_2$ is a leaf and $u_1$ is adjacent only to $u'$ in $A$, the only vertex which can separate $u_1$ and $u_2$ is $u'$. Therefore, $u'\in C'$. Thus, $C'\cap\{u,v,v',u'\}=\{u',v\}$ and $\{u,v'\}$ are not separated in $G$. Furthermore, we have $u''\not\in C'$ (otherwise $\{u,v'\}$ would be separated by $C'$ in $G$). Notice that since $C'$ separates $u'$ and $u$ in $G'$, we have $u_1$ in $C'$ and since $C'$ separates $u'$ and $u_1$, we have $u_2\in C'$. However, now the set $C=\{u\} \cup (C' \setminus \{u_1\})$ is an identifying code of claimed cardinality in $G$. Indeed, vertices $u',u,v\in C$ are adjacent and are thus separated by Lemma~\ref{LemP3}. Moreover, also vertices $u_1$ and $u_2$ have unique code neighborhoods. This contradicts the properties of $G$ and completes the proof of Claim~\ref{c:claimF}.~\smallqed

\medskip

Recall that $$1 \le |A^*| \le k \le b^* \le b.$$ By Claim~\ref{c:claimF}, we may assume that $A^* \subset A$, and so $1 \le |A^*| < |A| = a \le 2(\Delta - 1)$.  Let $\overline{A^*} = A \setminus A^*$. Thus, $|\overline{A^*}| = |A| - |A^*| \ge 1$. We note that either $|\overline{A^*}| \le \Delta - 2$ or $|\overline{A^*}| \ge \Delta - 1$. Further, either $\overline{A^*}$ contains a neighbor of $u$ and a neighbor of $v$ or $\overline{A^*} \subset N_{G^*}(u)$ or $\overline{A^*} \subset N_{G^*}(v)$. We proceed further with three claims.

\medskip
\begin{subclaim}
\label{c:claimG}
If $|\overline{A^*}| \le \Delta - 2$,  then there exists an identifying code in $G^*$ containing vertices $u$ and $v$ with cardinality at most
\[
\left( \frac{\Delta - 1}{\Delta}\right) n^*.
\]
\end{subclaim}
\proof Suppose that $|\overline{A^*}| \le \Delta - 2$. As observed earlier, $|A^*| \le b^*$. Thus in this case, $a = |A| = |A^*| + |\overline{A^*}| \le \Delta + b^* - 2$. By Claim~\ref{c:claimE}(b), we have $b - (\Delta - 1)b^* \le 0$. Let $x$ be an arbitrary vertex in $\overline{A^*}$ and let
\[
C^* = V(G^*) \setminus (B^* \cup \{x\}).
\]

The code $C^*$ is an identifying code of $G^*$. Indeed, $C^*$ is a dominating set in $G^*$. Furthermore, it is also connected and hence, by Lemma~\ref{LemP3}, separates every vertex in $C^*$. Furthermore, $x$ is the only vertex in $A$ which does not belong to set $C^*$ and it is separated by $\{u,v\}$ from vertices in $B$. Finally, vertices in $B^*$ are separated by $A^*\cup\{B\setminus B^*\}$. 
Next, we consider the cardinality of $C^*$:

\[
\begin{array}{lcl}
 |C^*| & = & (a + b + 2) - (b^*+1) \2 \\
& = & \displaystyle{ \left( \frac{\Delta - 1}{\Delta} \right) (a + b + 2) + \frac{1}{\Delta}(a + b - (\Delta - 2) - \Delta b^*) } \2 \\
& \le & \displaystyle{ \left( \frac{\Delta - 1}{\Delta} \right)n^* + \frac{1}{\Delta}\big((\Delta + b^* - 2) + b - \Delta + 2 - \Delta b^* \big) } \2 \\
& \le & \displaystyle{ \left( \frac{\Delta - 1}{\Delta} \right)n^* + \frac{1}{\Delta}\big(b - (\Delta - 1) b^* \big) } \2 \\
& \le & \displaystyle{ \left( \frac{\Delta - 1}{\Delta} \right)n^* }, \1 \\
\end{array}
\]
yielding the desired upper bound.~\smallqed

\medskip

\begin{subclaim}
\label{c:claimH}
If $\overline{A^*}$ contains a neighbor of $u$ and a neighbor of $v$,  then there exists an identifying code in $G^*$ containing vertices $u$ and $v$ with cardinality at most
\[
\left( \frac{\Delta - 1}{\Delta}\right) n^*.
\]
\end{subclaim}
\proof Suppose that $\overline{A^*}$ contains a neighbor, $u_1$, of $u$ and a neighbor, $v_1$, of $v$. By our earlier observations, $|\overline{A^*}| = a - |A^*| \le 2(\Delta - 1) - 1 = 2\Delta - 3$. Hence, $a = |A^*| + |\overline{A^*}| \le b^* + 2\Delta - 3$. By Claim~\ref{c:claimE}(b), we have $b - (\Delta - 1)b^* \le 0$. We now let
\[
C^* = V(G^*) \setminus (B^* \cup \{u_1,v_1\}).
\]

The set $C^*$ is an identifying code of $G^*$. Indeed, $C^*$ is a dominating set in $G^*$. Furthermore, it is also connected and hence, by Lemma~\ref{LemP3}, separates every vertex in $C^*$. Furthermore, $u_1$ and $v_1$ are the only vertices in $A$ which do not belong to set $C^*$ and are separated by $\{u,v\}$ from vertices in $B$ and from each other. Finally, vertices in $B^*$ are separated by $A^*\cup\{B\setminus B^*\}$. 
Next, we consider the cardinality of $C^*$:
\[
\begin{array}{lcl}
 |C^*| & = & (a + b + 2) - (b^*+2) \2 \\
& = & \displaystyle{ \left( \frac{\Delta - 1}{\Delta} \right) (a + b + 2) + \frac{1}{\Delta}(a + b - 2(\Delta - 1) - \Delta b^*) } \2 \\
& \le & \displaystyle{ \left( \frac{\Delta - 1}{\Delta} \right)n^* + \frac{1}{\Delta}\big((b^* + 2\Delta - 3) + b - 2(\Delta - 1) - \Delta b^* \big) } \2 \\
& \le & \displaystyle{ \left( \frac{\Delta - 1}{\Delta} \right)n^* + \frac{1}{\Delta}\big(b - (\Delta - 1) b^* - 1 \big) } \2 \\
& \le & \displaystyle{ \left( \frac{\Delta - 1}{\Delta} \right)n^* - \frac{1}{\Delta} } \2 \\
& < & \displaystyle{ \left( \frac{\Delta - 1}{\Delta} \right)n^* }, \1 \\
\end{array}
\]
yielding the desired upper bound.~\smallqed

\medskip

By Claim~\ref{c:claimG}, we may assume that $|\overline{A^*}| \ge \Delta - 1$, for otherwise there exists an identifying code of size at most $\left(\frac{\Delta-1}{\Delta}\right)n^*$ in $G^*$ containing vertices $u$ and $v$.

\medskip
\begin{subclaim}
\label{c:claimI}
If $\overline{A^*} \subseteq N_u$ or $\overline{A^*} \subseteq N_v$,   then there exists an identifying code in $G^*$ containing vertices $u$ and $v$ with cardinality at most
\[
\left( \frac{\Delta - 1}{\Delta}\right) n^*.
\]
\end{subclaim}
\proof Suppose that either $\overline{A^*} \subseteq N_u$ or $\overline{A^*} \subseteq N_v$. Renaming $u$ and $v$ if necessary, we may assume that $\overline{A^*} \subseteq N_v$, and so $\overline{A^*}$ contains no neighbor of $u$. Since $\overline{A^*} \subseteq N_v$, we have $|\overline{A^*}| \le |N_{G^*}(v)| - 1 \le \Delta - 1$. However by our earlier assumption due to Claim~\ref{c:claimG}, $|\overline{A^*}| \ge \Delta - 1$.  Therefore, $|\overline{A^*}| = \Delta - 1$, that is, $\deg_{G^*}(v) = \Delta$ and $\overline{A^*} = N_{G^*}(v) \setminus \{u\}$. Thus,
\[
a = |A| = |A^*| + |\overline{A^*}| \le b^* + \Delta - 1.
\]

Recall that $u'v'$ is an edge in $G^*$, where $u'$ is a neighbor of $u$ in $G^*$ different from $v$ and $v'$ is a neighbor of $v$ in $G^*$ different from $u$. Since $u' \in A^*$, the vertex $u'$ is adjacent to at least one vertex in $B$, and therefore $u'$ is adjacent to at most~$\Delta - 2$ vertices in $\overline{A^*}$. Let $v''$ be a vertex in $\overline{A^*}$ that is not adjacent to the vertex~$u'$. By Claim~\ref{c:claimE}(b), we have $b - (\Delta - 1)b^* \le 0$. We now let
\[
C^* = V(G^*) \setminus (B^* \cup \{v',v''\}).
\]

The set $C^*$ is an identifying code of $G^*$. Indeed, $C^*$ is a dominating set in $G^*$. Furthermore, it is also connected and hence, by Lemma~\ref{LemP3}, separates every vertex in $C^*$. Furthermore, $v'$ and $v''$ are the only vertices in $A$ which do not belong to set $C^*$ and are separated by $v$ from vertices in $B$ and by $u'$ from each other. Finally, vertices in $B^*$ are separated by $A^*\cup\{B\setminus B^*\}$. 
Next, we consider the cardinality of $C^*$:\[
\begin{array}{lcl}
 |C^*| & = & (a + b + 2) - (b^*+2) \2 \\
& = & \displaystyle{ \left( \frac{\Delta - 1}{\Delta} \right) (a + b + 2) + \frac{1}{\Delta}(a + b - 2(\Delta - 1) - \Delta b^*) } \2 \\
& \le & \displaystyle{ \left( \frac{\Delta - 1}{\Delta} \right)n^* + \frac{1}{\Delta}\big((b^* + \Delta - 1) + b - 2(\Delta - 1) - \Delta b^* \big) } \2 \\
& \le & \displaystyle{ \left( \frac{\Delta - 1}{\Delta} \right)n^* + \frac{1}{\Delta}\big(b - (\Delta - 1) b^* - (\Delta - 1)  \big) } \2 \\
& \le & \displaystyle{ \left( \frac{\Delta - 1}{\Delta} \right)n^* - \left( \frac{\Delta - 1}{\Delta} \right) } \2 \\
& < & \displaystyle{ \left( \frac{\Delta - 1}{\Delta} \right)n^* }, \1 \\
\end{array}
\]
yielding the desired upper bound.~\smallqed

\medskip

By Claims~\ref{c:claimF}, \ref{c:claimG}, \ref{c:claimH} and \ref{c:claimI}, we have shown that graph $G^*$ admits an identifying code of cardinality at most $\left(\frac{\Delta-1}{\Delta}\right)n^*$ containing vertices $u$ and $v$. Hence, Claim~\ref{claimGStar} follows.
~\smallqed

If $G^*=G$, then the theorem statement follows from Claim~\ref{claimGStar}. If $G^*\neq G$, then graph $G_{uv}$ contains a component of cardinality at least~3. Furthermore, by Claim~\ref{c:claimC}, graph $G_{uv}$ does not contain any components that belong to the set $\mathcal{F}_\Delta$.  Hence, we may assume that $G_{uv}$ contains a component of order at least~$3$. Next, we finalize the proof by showing that also in this case, $\ID(G)\le \left( \frac{\Delta - 1}{\Delta}\right) n$.

Let us denote by $\mathcal{K}$ the set of components of order at least~$3$ in $G_{uv}$. By Claim~\ref{c:claimC}, $K \notin \mathcal{F}_{\Delta}$ for any $K\in \mathcal{K}$.
Applying the inductive hypothesis to a component $K\in \mathcal{K}$ of maximum degree $\Delta_K$, component $K$ satisfies
\begin{equation}
\label{Eq1F}
\ID(K) \le \left( \frac{\Delta_K - 1}{\Delta_K}\right) n(K) \le \left( \frac{\Delta - 1}{\Delta}\right) n(K).
\end{equation} Let $C_K$ be an identifying code in $K$ with minimum cardinality.

Observe that $G^* = G - \bigcup_{K\in \mathcal{K}} V(K)$. Let $C^*$ be an identifying code of $G^*$ with cardinality at most $\left(\frac{\Delta-1}{\Delta}\right)n^*$ containing vertices $u$ and $v$, that exists by Claim~\ref{claimGStar}.

Let us next consider set $C=C^*\cup\left(\bigcup_{K\in \mathcal{K}}C_K\right)$. 
Set $C$ is an identifying code in $G$ since set $C^*$ contains vertices $u$ and $v$. Indeed, since $C$ a union of multiple identifying codes, each identifying code dominates and pairwise separates the vertices within the corresponding component. Thus, the only problems might be between vertices of different components. In particular, in the case the component $K$ contains vertex $z\in C_K$ such that $N[z]\cap C_K=\{z\}$, component $G^*$ contains vertex $y\in C^*$ such that $N[y]\cap C^*=\{y\}$, and $y$ and $z$ are adjacent in $G$. However, this is not possible since $y\in A$ and thus, $N[y]\cap C^*\cap\{u,v\}\neq\emptyset$.

Notice that we have \[
\ID(G) \le |C|= |C^*| +\sum_{K\in\mathcal{K}}|C_K| \le \left( \frac{\Delta - 1}{\Delta}\right) n^* + \sum_{K\in\mathcal{K}}\left( \frac{\Delta - 1}{\Delta}\right) n(K) = \left( \frac{\Delta - 1}{\Delta}\right)n.
\]

This completes the proof of Theorem~\ref{thm:trianglefree}.~\hfill\QED

\medskip
As an immediate consequence of Theorem~\ref{thm:trianglefree}, we have the following corollary.

\begin{corollary}
\label{cor:trianglefree1}
If $G$ is a connected, identifiable, triangle-free graph of order $n$, then the following holds.  \\[-18pt]
\begin{enumerate}
\item[{\rm (a)}] If $G$ is cubic, then  $\ID(G) \le \frac{2}{3}n$.  \1
\item[{\rm (b)}] If $G$ is subcubic and $n \ge 23$, then $\ID(G) \le \frac{2}{3}n$. \1
\item[{\rm (c)}] If $G$ has maximum degree~$\Delta$ where $\Delta \ge 4$ is fixed and $n \ge \Delta + 2$, then $\ID(G) \le \left( \frac{\Delta - 1}{\Delta}\right) n$.
\end{enumerate}
\end{corollary}

\section{Beyond triangle-free graphs}\label{sec:beyond}

We now apply Theorem~\ref{thm:trianglefree} to graphs having triangles and obtain a bound weaker than the conjectured one, as follows.

\begin{corollary}\label{cortTriangles}
If $G$ is a connected identifiable graph of order $n \ge 3$ with maximum degree $\Delta\ge 3$ such that $G$ can be made triangle-free by deleting $t$ edges, then
\[
\ID(G)\le \left( \frac{\Delta-1}{\Delta} \right) n+4t+\frac{1}{\Delta}.
\]
\end{corollary}
\begin{proof}
We prove the claim by induction on $t$. Let us assume that $E_t\subseteq E(G)$, with $|E_t|=t$, is a smallest set of edges which we can remove from $G$ so that $G_t=G-E_t$ is triangle-free. Observe that $G_t$ is connected and identifiable, since the only connected triangle-free graph that is not identifiable is $K_2$, but we assume here that $n\ge 3$. By Table~\ref{TableConstVal} and Theorem~\ref{thm:trianglefree}, the claim holds when $G$ is triangle-free, that is, for $t=0$. Assume next that the claim holds for every $t\le t'$ and let $t=t'+1$. Assume that $G_t$ has maximum degree $\Delta_t\le \Delta$. We have $\Delta_t\ge 2$ since $G_t$ is connected and $n\ge 3$.

Let $C_t$ be an optimal identifying code in $G_t$. By Theorem~\ref{thm:trianglefree}, if $\Delta_t\ge 3$, we have
\[
|C_t|\le \left( \frac{\Delta_t-1}{\Delta_t} \right) n + \frac{1}{\Delta_t}\le \left( \frac{\Delta-1}{\Delta} \right) n + \frac{1}{\Delta}.
\]
If $\Delta_t=2$, by Corollary~\ref{cor_paths & cycles}(e)-(f) and $n\notin\{4,7\}$, then $|C_t| \le \frac{2}{3}n \le \left( \frac{\Delta-1}{\Delta} \right)n$ since $\Delta\ge 3$.
If $n=4$, $|C_t|=3\le \left( \frac{\Delta-1}{\Delta} \right) n + \frac{1}{\Delta}$ and if $n=7$, $|C_t|\le 5\le \left( \frac{\Delta-1}{\Delta} \right) n+\frac{1}{\Delta}$ (again since $\Delta\ge 3$). Let edge $uv\in E_t$ and let us consider graph $G_t$ together with edge $uv$, denoted by $G_t'$. Observe that if $C_t$ is not an identifying code in $G_t'$, then the addition of edge $uv$ modified some code-neighborhoods (this implies that $u$ or $v$ is in $C_t$). Therefore, there are at most four vertices which are no longer separated (possibly, vertex $u$ together with some vertex $u'$, and possibly, vertex $v$ with some vertex $v'$). As we add these edges back one at a time, each time we create at most four new vertices among which some vertex-pairs are not separated by $C_t$. Thus, in the graph $G$, there is a set $S$ of at most $4t$ vertices in which some vertex-pairs are not separated by $C_t$. However, since $G$ is identifiable, $S$ is $V(G)$-identifiable and thus, by Lemma~\ref{lemXYID}, we can find an $(S,V(G))$-identifying code $C_S$ of size at most $|S|\le 4t$. The set $C_t\cup C_S$ is an identifying code of $G$ of the desired size, proving the claim.
\end{proof}

\section{Conclusion}\label{sec:conclu}
We proved Conjecture~\ref{conj_G Delta_UB} for all triangle-free graphs. In fact, we proved (see Theorem~\ref{thm:trianglefree}) the following stronger result. Let $G$ be a connected, identifiable, triangle-free graph of order $n \ge 3$ with maximum degree $\Delta$. If $\Delta \ge 4$, then we proved that $\ID(G) \le \left( \frac{\Delta - 1}{\Delta}\right) n$, except for one exceptional graph, namely the star $K_{1,\Delta}$. Moreover, if $\Delta \le 3$, then $\ID(G) \le \frac{2}{3}n$, unless $G$ belongs to a forbidden family that contains fifteen graphs (all of order at most~$22$): $P_4$, $C_4$, $C_7$ and the twelve trees of maximum degree~3 from $\mathcal{T}_3$.

In the special case when $G$ is a triangle-free cubic graph, this implies that $\ID(G) \le \frac{2}{3}n$ always holds. This establishes a best possible upper bound for triangle-free cubic graphs, since $\ID(K_{3,3})=4$ (we do not know other cubic graphs for which this bound is tight exist).
The previously best known upper bounds (prior to this paper) when $G$ is triangle-free subcubic and cubic were, $\ID(G) \le \frac{8}{9}n$ and $\ID(G) \le \frac{5}{6}n$, respectively (see~\cite[Corollary 4.46]{FlorentPhD}; the proof used the technique developed in~\cite{foucaud2012size}).

Towards a positive resolution of  Conjecture~\ref{conj_G Delta_UB}, it would be interesting to prove it for all cubic graphs.

When it comes to general graphs, the list of exceptional graphs is larger than $\mathcal{F}_\Delta$. Indeed, as mentioned in Table~\ref{TableConstVal}, the complements of half-graphs and related constructions defined in~\cite{FGKNPV11} do require $c>0$. 
Nevertheless, those constructions have maximum degree $\Delta$ very close to the number $n$ of vertices ($n-1$ or $n-2$), and, for any given $\Delta$, there is only a finite number of such examples. Thus, it is possible that even for the general case, if the conjecture is true, the list of graphs requiring $c>0$ is also finite for every fixed value of $\Delta$.

As seen in Table~\ref{TableConstVal}, all graphs known to us that require $c>0$ have $c\le 3/2$ (which is reached only by odd cycles), and when $\Delta\ge 3$, in fact $c\le 1/3$. Are there graphs that require higher values of $c$? Another way to formulate these constants is in terms of the maximum degree. Do there exist any graphs that require the constant $c$ to be larger than $3/\Delta$? By our results, such graphs would necessarily contain triangles. Note that it seems necessary to understand those graphs needing $c>0$, in order to prove the conjecture.

\bibliographystyle{abbrv}
\bibliography{Bibliography}

\end{document}